\documentclass[11pt]{amsart}
\usepackage[a4paper,hmargin=1.9cm,top=2.3cm,bottom=2.3cm]{geometry}
\usepackage{mathpazo}
\usepackage{amsmath,amssymb,amsthm,mathtools}
\usepackage{xcolor}
\usepackage{microtype}
\usepackage[colorlinks,linkcolor=blue!60!black,citecolor=green!45!black,urlcolor=blue!60!black]{hyperref}

\newif\ifhighlight
\highlightfalse
\newcommand{\chg}{\ifhighlight\color{blue}\fi}

\newcommand{\cO}{\mathcal{O}}
\newcommand{\vp}{\varpi}
\newcommand{\one}{\mathbb{1}}
\newcommand{\C}{\mathbb{C}}
\newcommand{\F}{\mathbb{F}}
\newcommand{\Fq}{\mathbb{F}_q}
\newcommand{\M}{\mathrm{M}}
\newcommand{\Oh}{\widehat{\cO_2^{\times}}}
\newcommand{\ct}{\widetilde{\chi}}
\newcommand{\Sym}{\mathfrak{S}}
\DeclareMathOperator{\GL}{GL}
\DeclareMathOperator{\Ind}{Ind}
\DeclareMathOperator{\Res}{Res}
\DeclareMathOperator{\Hom}{Hom}
\DeclareMathOperator{\End}{End}
\DeclareMathOperator{\tr}{tr}
\DeclareMathOperator{\diag}{diag}
\DeclareMathOperator{\St}{St}

\theoremstyle{plain}
\newtheorem{theorem}{Theorem}[section]
\newtheorem{lemma}[theorem]{Lemma}
\newtheorem{proposition}[theorem]{Proposition}
\newtheorem{corollary}[theorem]{Corollary}
\newtheorem*{theoremA}{Theorem A}
\newtheorem*{theoremB}{Theorem B}
\theoremstyle{definition}
\newtheorem{definition}[theorem]{Definition}
\theoremstyle{remark}
\newtheorem{remark}[theorem]{Remark}
\numberwithin{equation}{section}
\numberwithin{table}{section}
\newtheorem{Conjecture}[theorem]{Conjecture}

\begin{document}

\title{Structure of principal series representations of $\GL_3(\cO_2)$}
\author{Ankita Parashar}
\address{Department of Mathematics, IIT Bombay, Powai, Maharashtra 400076, India}
\email{ankitaparashar216@gmail.com}

\begin{abstract}
Let $F$ be a non-Archimedean local field with ring of integers $\cO$ and uniformizer $\vp$, and for a positive integer $\ell$ let $\cO_\ell:=\cO/(\vp^\ell)$.
For a character $\chi$ of the diagonal torus, let $I(\chi)=\Ind_B^{\GL_n(\cO_\ell)}\chi$ denote the principal series representation, where $B$ is the Borel subgroup of upper triangular matrices.
For $\GL_2$, the decomposition of $I(\chi)$ is classical when $\ell=1$, and was obtained for all $\ell$ in work of Casselman.
In this paper we give a complete decomposition of $I(\chi)$ for $\GL_3(\cO_2)$ and every character $\chi$: we determine its irreducible constituents, their multiplicities.
Unlike for $\GL_2$, $I(\chi)$ need not be multiplicity free.
Nevertheless, we show that $I(\chi)$ always contains exactly one regular constituent and it occurs with multiplicity one.
We conjecture that the same holds in general: for all $n,\ell\ge2$ and every character $\chi$, the principal series representation $I(\chi)$ of $\GL_n(\cO_\ell)$ contains exactly one irreducible regular representation, and it occurs with multiplicity one.
\end{abstract}

\maketitle

\section{Introduction}

Let $F$ be a finite extension of $\mathbb{Q}_p$, the field of $p$-adic numbers or $\F_p((t))$, the field of formal Laurent series over the finite field $\F_p$. Let $\cO$ be its ring of integers, $\vp$ a uniformizer and $\Fq=\cO/(\vp)$ the residue field. For $\ell\ge1$ the quotient $\cO_\ell:=\cO/(\vp^\ell)$ is a finite principal ideal local ring of length $\ell$; for example $\cO_\ell=\mathbb{Z}/p^\ell\mathbb{Z}$ or $\F_p[t]/(t^\ell)$.

The group $\GL_n(\cO)$ is a maximal compact subgroup of $\GL_n(F)$, and every irreducible smooth representation of $\GL_n(\cO)$ factors through $\GL_n(\cO_\ell)$ for some $\ell$. The representation theory of the finite groups $\GL_n(\cO_\ell)$ is therefore closely tied to that of $\GL_n(F)$. For $\ell=1$ it is the classical theory of $\GL_n(\Fq)$. For $\ell=2$, the classification of all irreducible representations of $\GL_n(\cO_\ell)$ is given in~\cite{Singla10} and recently in~\cite{OPS25}, the complete description of irreducible representations of $\GL_3(\cO_\ell)$ is known for all $\ell$. For $\ell>2$ it is much less understood, and a complete description of the irreducible representations is not known in full generality. Important contributions include the work of Hill~\cite{Hill95}, Krakovski-Onn-Singla~\cite{KOS18}, Singla~\cite{Singla10}, Stasinski-Stevens~\cite{SS17} and Stasinski~\cite{Stasinski17}; see also the work of Leigh-Cliff-Wen~\cite{BCW10} on the character values of $\GL_2(\mathbb{Z}/p^\ell\mathbb{Z})$.

Let $B=TU$ be the subgroup of upper triangular matrices in $\GL_n(\cO_\ell)$, where $T$ is the diagonal torus and $U$ the subgroup of upper unitriangular matrices. A character $\chi=\chi_1\otimes\cdots\otimes\chi_n$ of $T$, with $\chi_i$ characters of $\cO_\ell^\times$, extends to $B$ through the projection $B\twoheadrightarrow B/U\cong T$, and we call $\Ind_B^{\GL_n(\cO_\ell)}\chi$ a \emph{principal series representation}. The basic problem is to determine its irreducible constituents and their multiplicities.

For $\ell=1$ the answer is classical and follows from the Bruhat decomposition. For $\ell\ge2$ principal series representations are typically highly reducible, and far fewer tools are available: analogues of the Bruhat decomposition are known only in special cases~\cite{OPV06}, and the character theory of $\GL_n(\cO_\ell)$ is not well understood. The problem is also connected with several others. For instance, multiplicity freeness of principal series representations is closely related to the question of whether $(\GL_n(\cO_\ell),B)$ is a Gelfand pair~\cite{GS25}. For $n=2$ the principal series of $\GL_2(\cO_\ell)$ is multiplicity free~\cite{Cas73}. For $n\ge3$, however, Dagar~\cite{Dag26} recently showed that, for every $\ell$, $\Ind_B^{\GL_n(\cO_\ell)}\chi$ is not multiplicity free when all the components of $\chi$ are equal; in particular, $(\GL_n(\cO_\ell),B)$ is not a Gelfand pair. For $\GL_3(\cO_\ell)$ he also gave an irreducibility criterion and decomposed the principal series under some conditions.

In this paper we treat every character $\chi$ of $\GL_3(\cO_2)$, including the cases not covered by~\cite{Dag26}. Our results complement the following works. Onn and Singla, Campbell and Nevins study the restriction of a principal series representation $\Ind_{B(F)}^{\GL_3(F)}\chi$ of $\GL_3(F)$ to its maximal compact subgroup $\GL_3(\cO)$. By the Iwasawa decomposition $\GL_3(F)=B(F)\GL_3(\cO)$, this restriction is $\Ind_{B(\cO)}^{\GL_3(\cO)}\chi$, an infinite-dimensional representation of the compact group $\GL_3(\cO)$. Its vectors fixed by $K^\ell=\ker\bigl(\GL_3(\cO)\to\GL_3(\cO_\ell)\bigr)$ form exactly our finite-level representation $\Ind_{B(\cO_\ell)}^{\GL_3(\cO_\ell)}\chi$, whenever $\chi$ is trivial on $T(1+\vp^\ell\cO)$. Thus a complete decomposition of the restriction gives the decomposition of $\Ind_{B(\cO_\ell)}^{\GL_3(\cO_\ell)}\chi$ for every $\ell$, and conversely. However, translating a description of the restriction into an explicit list of constituents at a fixed level $\ell$ from these two works is not automatic.

Onn and Singla~\cite{OS14} give a complete description of the irreducible constituents of the restriction when $\chi$ is unramified, that is, trivial on $T(\cO)$. In our setting this covers $\Ind_{B(\cO_\ell)}^{\GL_3(\cO_\ell)}\one$ for every $\ell$, and hence, after twisting by a character of the determinant, every $\chi$ with $\chi_1=\chi_2=\chi_3$. At level $2$ this is the single case $\Ind_{B(\cO_2)}^{\GL_3(\cO_2)}\one$, up to twist.

Campbell and Nevins study the same restriction for unramified characters~\cite{CN09} and for ramified characters~\cite{CN10}. In the ramified case, assuming $p>2$ and $q>3$, they write the restriction as a direct sum of quotients of representations induced from subgroups defined by filtrations of the Iwahori subgroup. The shape of this decomposition depends on the ramification of the characters $\chi_i\chi_j^{-1}$. They prove that many of these quotients are irreducible \cite[Theorems~4.1--4.3]{CN10}, but not all of them are \cite[Example~5.1]{CN10}. So a complete list of irreducible constituents with their multiplicities is not obtained, nor are explicit models given for the constituents of the reducible quotients. In contrast, for $\GL_3(\cO_2)$ we treat every character $\chi$, ramified or not, and for each $\chi$ we (i) determine all irreducible constituents of $I(\chi)$ together with their dimensions and multiplicities; (ii) realise each constituent explicitly as a representation induced from a character of an explicit subgroup; (iii) compute $\dim\End(I(\chi))\in\{1,2,3,7,10,18\}$; and (iv) show that $I(\chi)$ has exactly one regular constituent, which occurs with multiplicity one.

In the work of Casselman, see \cite[Proposition~1]{Cas73}, we get the decomposition of principal series representations of $\GL_2(\cO_\ell)$. The following theorem summarizes work of Casselman for $\GL_2(\cO_\ell)$. In the statement below we use the convention $1+\vp^0\cO_\ell:=\cO_\ell^\times$.

\begin{theorem}[\cite{Cas73}]\label{thm:Casselman}
Let $\chi_1,\chi_2$ be characters of $\cO_\ell^\times$ and $\chi=\chi_1\otimes\chi_2$.
\begin{enumerate}
\item[(1)] The representation $\Ind_B^{\GL_2(\cO_\ell)}\chi$ has $\ell+1$ irreducible constituents if and only if $\chi_1=\chi_2$.
\item[(2)] Suppose that $\chi_1\neq\chi_2$, and let $0\le k\le\ell-1$ be the unique integer such that $\chi_1$ and $\chi_2$ agree on $1+\vp^{k+1}\cO_\ell$ but not on $1+\vp^k\cO_\ell$. Then
\[
\Ind_B^{\GL_2(\cO_\ell)}\chi\cong\bigoplus_{i=k}^{\ell-1}\sigma_i,
\]
where the $\sigma_i$ are irreducible and pairwise non-isomorphic, $\dim\sigma_k=q^k(q+1)$ and $\dim\sigma_i=q^{i-1}(q^2-1)$ for $k+1\le i\le\ell-1$.
\end{enumerate}
\end{theorem}

In particular, the principal series representations of $\GL_2(\cO_\ell)$ are multiplicity free, and each of them contains exactly one irreducible regular representation (see Definition \ref{def:regular}), with multiplicity one.
Consider the principal series representation
\[
I(\chi):=\Ind_B^{\GL_3(\cO_2)}\chi,\qquad\chi=\chi_1\otimes\chi_2\otimes\chi_3,
\]
of $\GL_3(\cO_2)$. For characters $\lambda,\mu$ of $\cO_2^\times$ we write $\lambda\sim\mu$ if $\lambda$ and $\mu$ agree on $1+\vp\cO_2$. The first step is an explicit description of the double coset space $B\backslash\GL_3(\cO_2)/B$, which has $18$ elements for every $q$, see Proposition~\ref{prop:18}. Combined with Mackey theory, this gives the dimension of the endomorphism algebra of $I(\chi)$, see Proposition~\ref{prop:Dchi}.

\begin{theoremA}
The dimension $D(\chi):=\dim\End_{\GL_3(\cO_2)}(I(\chi))$ is given by
\[
D(\chi)=\begin{cases}
1 & \text{if } \chi_i\not\sim\chi_j \text{ for all } i\neq j,\\
2 & \text{if } \chi_i\sim\chi_j\not\sim\chi_k \text{ and } \chi_i\neq\chi_j,\\
3 & \text{if } \chi_i=\chi_j\not\sim\chi_k,\\
7 & \text{if } \chi_1\sim\chi_2\sim\chi_3 \text{ and the } \chi_i \text{ are pairwise distinct},\\
10 & \text{if } \chi_1\sim\chi_2\sim\chi_3 \text{ and } \chi_i=\chi_j\neq\chi_k,\\
18 & \text{if } \chi_1=\chi_2=\chi_3,
\end{cases}
\]
where $\{i,j,k\}=\{1,2,3\}$.
\end{theoremA}

In each of these cases we then decompose $I(\chi)$ explicitly into irreducible representations (Sections~\ref{sec:D23}--\ref{sec:D1018}). Moreover, we get the following corollary.

\begin{corollary}\label{cor:mf}
$I(\chi)$ is multiplicity free when $D(\chi)\le7$, and it is not multiplicity free when $D(\chi)\in\{10,18\}$.
\end{corollary}

To state our second main result, let $K=\ker\bigl(\GL_3(\cO_2)\to\GL_3(\Fq)\bigr)$. This is an abelian normal subgroup of $\GL_3(\cO_2)$ isomorphic to $\M_3(\Fq)$. Fix a character $\psi$ of $\cO_2$ which is non-trivial on $\vp\cO_2$. Then every character of $K$ has the form $\psi_\beta\colon I_3+\vp Y\mapsto\psi(\vp\tr(\beta Y))$ for a unique $\beta\in\M_3(\Fq)$. An irreducible representation of $\GL_3(\cO_2)$ is called \emph{regular} if it lies over a character $\psi_\beta$ with $\beta$ regular, that is, with minimal polynomial of $\beta$ is equal to its characteristic polynomial of $\beta$.

\begin{theoremB}
For every $\chi$, the representation $I(\chi)$ contains exactly one irreducible regular representation, and it occurs in $I(\chi)$ with multiplicity one.
\end{theoremB}

This is the analogue of the corresponding property of the principal series of $\GL_2(\cO_\ell)$ recalled above. Its dimension is determined by which of the $\chi_i$ are related by $\sim$, see Corollary~\ref{cor:regdim} for more details.
Theorem~\ref{thm:Casselman} and Theorem B suggest that this property holds for all $n$ and $\ell$.

\begin{Conjecture}
Let $n\ge2$ and $\ell\ge2$. For every character $\chi$ of $T$, the principal series representation $\Ind_B^{\GL_n(\cO_\ell)}\chi$ contains exactly one irreducible regular representation, and it occurs with multiplicity one.
\end{Conjecture}

The key idea of the proofs is to use suitable chains of subgroups between $B$ and $\GL_3(\cO_2)$, to apply transitivity of induction, and to describe the representations that occur at each step.

For $\ell>2$, the number of double cosets in $B\backslash\GL_3(\cO_\ell)/B$ is known to be a polynomial in $q$~\cite{OPV06}. For $\ell=2$, in contrast, the double coset structure is independent of $q$, which makes this case particularly amenable to a detailed analysis. Even when the double cosets are known, the decomposition of $\Ind_B^{\GL_3(\cO_\ell)}\chi$ remains unknown in general. The present paper gives a complete description for $\GL_3(\cO_2)$ and highlights several phenomena that are expected to play a role for larger groups.

\subsection*{Organization}
Section~\ref{sec:prelim} fixes notation and state some key results that are crucial for the entire article.
In Section~\ref{sec:dc} we determine $B\backslash\GL_3(\cO_2)/B$ and compute $D(\chi)$. The cases $D(\chi)\in\{2,3\}$ are treated in Section~\ref{sec:D23}, the case $D(\chi)=7$ in Section~\ref{sec:D7}, and the cases $D(\chi)\in\{10,18\}$ in Section~\ref{sec:D1018}. 
The regular constituent of $I(\chi)$ is determined in Section~\ref{subsec:regular}.

\section{Preliminaries}\label{sec:prelim}

\subsection{Notation}
For integers $1\le m<\ell$, the natural projection $\cO_\ell\twoheadrightarrow\cO_m$ induces a surjection $\GL_3(\cO_\ell)\twoheadrightarrow\GL_3(\cO_m)$. The $m$-th principal congruence subgroup of $\GL_3(\cO_\ell)$ is
\[
K_\ell^m:=\ker\bigl(\GL_3(\cO_\ell)\twoheadrightarrow\GL_3(\cO_m)\bigr),
\]
and we obtain the filtration $\{I\}\subseteq K_\ell^{\ell-1}\subseteq\cdots\subseteq K_\ell^1\subseteq\GL_3(\cO_\ell)$. For $1\le i\le\ell$, put $B_i:=B\cdot K_\ell^i$. This gives the filtration of subgroups
\[
B=B_\ell\subseteq B_{\ell-1}\subseteq\cdots\subseteq B_1\subseteq B_0:=\GL_3(\cO_\ell).
\]
From now on we take $\ell=2$ and use the following notation.
\begin{enumerate}
\item[(1)] $G:=\GL_3(\cO_2)$, $K:=K_2^1=I_3+\vp\M_3(\cO_2)$ and $B_1:=B\cdot K$. It is easy to see that being an abelian group, $K$ can be identified with the additive group $\M_3(\Fq)$. The subgroup $K$ is normal in $G$ and $G/K\cong\GL_3(\Fq)$.
\item[(2)] For a subgroup $H\le G$ we write $\overline{H}$ for its image in $\GL_3(\Fq)$.
\item[(3)] $I_3$ is the identity matrix and $E_{ij}$ is the matrix with $1$ in the $(i,j)$-entry and $0$ elsewhere.
\item[(4)] For $w\in\Sym_3$, let $P_w$ be the permutation matrix with $P_we_j=e_{w(j)}$. Then $P_w^{-1}E_{ij}P_w=E_{w^{-1}(i)\,w^{-1}(j)}$.
\item[(5)] $\Oh$ denotes the group of characters of $\cO_2^\times$. For $\lambda,\mu\in\Oh$ we write $\lambda\sim\mu$ if $\lambda|_{1+\vp\cO_2}=\mu|_{1+\vp\cO_2}$.
\item[(6)] For subgroups $H,K$ of a finite group $G$, we write $[K\backslash G/H]$ for a set of representatives of the double cosets $K\backslash G/H$.
\end{enumerate}
Since $\cO_2^\times/(1+\vp\cO_2)\cong\Fq^\times$, two characters $\lambda,\mu$ satisfy $\lambda\sim\mu$ if and only if $\lambda\mu^{-1}$ factors through $\cO_2^\times\twoheadrightarrow\Fq^\times$.

\subsection{Mackey theory}
Let $H$ and $K$ be subgroups of a finite group $G$ and let $\rho$ be a representation of $H$. For $s\in G$ put $H_s:=sHs^{-1}\cap K$ and define the representation $\rho^s$ of $H_s$ by
\begin{equation}\label{eq:conj}
\rho^s(x)=\rho(s^{-1}xs)\qquad(x\in H_s).
\end{equation}
Mackey theory (see~\cite{Serre}) states that
\begin{equation}\label{eq:mackey}
\Res^G_K\Ind_H^G\rho\;\cong\bigoplus_{s\in[K\backslash G/H]}\Ind_{H_s}^K\rho^s.
\end{equation}
Combined with Frobenius reciprocity this gives, for characters $\lambda$ of $H$ and $\mu$ of $K$,
\begin{equation}\label{eq:hom}
\Hom_\Gamma\bigl(\Ind_H^G\lambda,\Ind_K^G\mu\bigr)\;\cong\bigoplus_{s\in[K\backslash G/H]}\Hom_{sHs^{-1}\cap K}\bigl(\lambda^s,\mu\bigr).
\end{equation}
Each summand has dimension at most one, and it is non-zero if and only if $\lambda^s=\mu$ on $sHs^{-1}\cap K$. In particular,
\begin{equation}\label{eq:end}
\End_G\bigl(I(\chi)\bigr)\;\cong\bigoplus_{\delta\in[B\backslash G/B]}\Hom_{\delta B\delta^{-1}\cap B}\bigl(\chi^\delta,\chi\bigr).
\end{equation}

\subsection{Regular representations of \texorpdfstring{$\GL_n(\cO_l)$}{GLn(Ol)}}

\begin{definition}\label{def:regular}
\begin{enumerate}
\item[(1)] Let $x\in M_n(\Fq)$. Then $x$ is called \emph{regular} if its characteristic polynomial is the same as its minimal polynomial.
\item[(2)] For $x\in M_n(\Fq)$, the character $\phi_x\colon M_n(\Fq)\to\C^\times$ is called \emph{regular} if $x$ is a \emph{regular} element.
\item[(3)] An irreducible representation of $\GL_n(\mathfrak{o}_l)$ is called \emph{regular} if its restriction to $K_l^{l-1}\cong M_n(\Fq)$ contains a character which is regular.
\end{enumerate}
\end{definition}

In the representation theory of $\GL_n(\mathfrak{o}_l)$, \emph{regular representations} play a central role. In what follows, we recall the construction of regular representations of $\GL_n(\mathfrak{o}_l)$ under the assumption that $l$ is even \cite[Theorem~4.1]{Hill95}.

\begin{theorem}\label{thm:Hill}
Let $l>1$ be an even integer and set $m=l/2$. Let $\rho\colon K_l^m\to\C^\times$ be a character such that $\rho|_{K_l^{l-1}}$ is a regular character. Let the inertia group of $\rho$ be
\[
I(\rho):=\left\{g\in\GL_n(\mathfrak{o}_l):\rho(g^{-1}xg)=\rho(x)\ \forall\,x\in K_l^m\right\}.
\]
Then the following hold.\\
(1) The character $\rho$ extends to its inertia group $I(\rho)$. If $\tilde\rho$ is an extension of $\rho$ to $I(\rho)$ then $\Ind_{I(\rho)}^{\GL_n(\mathfrak{o}_l)}\tilde\rho$ is an irreducible regular representation of $\GL_n(\mathfrak{o}_l)$.\\
(2) For a given regular representation $\pi$ of $\GL_n(\mathfrak{o}_l)$ such that $\rho$ appears in $\pi|_{K_l^m}$, there exists a character $\tilde\rho$ of $I(\rho)$ extending $\rho$ such that $\pi=\Ind_{I(\rho)}^{\GL_n(\mathfrak{o}_l)}\tilde\rho$.
\end{theorem}

\begin{remark}
For $l=2m$, $K_l^m\cong M_n(\cO_m)$ and any character of $M_n(\cO_m)$ is $\phi_B$ for some $B\in M_n(\cO_m)$. For a regular representation $\pi$ of $\GL_n(\cO_l)$, if $\phi_B$ appears in $\pi|_{K_l^m}$ then $\phi_{B'}$ appears in $\pi|_{K_l^m}$ if and only if $B'=gBg^{-1}$ for some $g\in\GL_n(\cO_m)$.
\end{remark}

\begingroup\chg
\subsection{A key lemma}\label{subsec:split}
We fix an additive character $\psi\colon\cO_\ell\to\C^\times$ with $\psi|_{\vp^{\ell-1}\cO_\ell}\neq\one$.
Note that $\psi(\vp^{\ell-1}y)$ depends only on $y\bmod\vp$.
We use the following lemma.

\begin{lemma}\label{lem:split}
Let $H\subseteq H^+$ be finite groups with $[H^+:H]=q$, and let $\lambda^+$ be a character of $H^+$.
Let $H'\trianglelefteq H^+$ be a normal subgroup with $H^+=HH'$ and $[H^+:H']=q-1$. Let $\theta$ be a non-trivial character of $H'$ with $\theta|_{H\cap H'}=\one$ and $\theta^x\neq\theta$ for every $x\in H^+\smallsetminus H'$.
Then
\[
\Ind_H^{H^+}\bigl(\lambda^+|_H\bigr)\cong\lambda^+\oplus\Ind_{H'}^{H^+}\bigl(\lambda^+|_{H'}\theta\bigr),
\]
and the second summand is irreducible of dimension $q-1$.
\end{lemma}

\begin{proof}
Since $H^+=H'H$, there is only one $(H',H)$-double coset in $H^+$, so by~\eqref{eq:mackey} we have $\Res^{H^+}_{H'}\Ind_H^{H^+}\one\cong\Ind_{H\cap H'}^{H'}\one$. As $\theta|_{H\cap H'}=\one$, Frobenius reciprocity shows that $\theta$ occurs in $\Ind_{H\cap H'}^{H'}\one$, and hence
\[
\Hom_{H^+}\bigl(\Ind_{H'}^{H^+}\theta,\ \Ind_H^{H^+}\one\bigr)\cong\Hom_{H'}\bigl(\theta,\ \Res^{H^+}_{H'}\Ind_H^{H^+}\one\bigr)\neq0.
\]
Since $H'$ is normal in $H^+$, Mackey theory gives $\End_{H^+}\bigl(\Ind_{H'}^{H^+}\theta\bigr)\cong\bigoplus_{x\in H^+/H'}\Hom_{H'}(\theta^x,\theta)$, and by assumption only the trivial coset contributes. Hence $\sigma:=\Ind_{H'}^{H^+}\theta$ is irreducible, of dimension $[H^+:H']=q-1$. It is not isomorphic to $\one$: for $q>2$ the dimensions differ, and for $q=2$ we have $H'=H^+$ and $\sigma=\theta\neq\one$. Since $\one$ also occurs in $\Ind_H^{H^+}\one$, the representation $\one\oplus\sigma$ embeds into $\Ind_H^{H^+}\one$, and as $1+(q-1)=q=\dim\Ind_H^{H^+}\one$ we get $\Ind_H^{H^+}\one\cong\one\oplus\sigma$. Tensoring with $\lambda^+$ and using $\lambda^+\otimes\Ind_H^{H^+}\one\cong\Ind_H^{H^+}(\lambda^+|_H)$ and $\lambda^+\otimes\Ind_{H'}^{H^+}\theta\cong\Ind_{H'}^{H^+}(\lambda^+|_{H'}\theta)$ gives the stated isomorphism. The second summand is the twist of the irreducible representation $\sigma$ by the character $\lambda^+$, hence it is irreducible.
\end{proof}

\section{The double cosets \texorpdfstring{$B\backslash G/B$}{B\textbackslash G/B} and the dimension of \texorpdfstring{$\End_G(I(\chi))$}{End(I(chi))}}\label{sec:dc}

In this section we compute $D(\chi)=\dim\End_G(I(\chi))$ using Mackey theory. To this end we first describe the double coset space $B\backslash G/B$ explicitly.

\begin{lemma}\label{lem:Bruhat}
The permutation matrices $\{P_w:w\in\Sym_3\}=\{I_3,P_{(12)},P_{(13)},P_{(23)},P_{(123)},P_{(132)}\}$ form a set of representatives for $\overline{B}\backslash\GL_3(\Fq)/\overline{B}$.
\end{lemma}

\begin{proof}
This is the Bruhat decomposition of $\GL_3(\Fq)$; see \cite[Exercise~4.1.18]{Bump}.
\end{proof}

\begin{proposition}\label{prop:18}
The set $B\backslash G/B$ has $18$ elements. A set of representatives is $\Gamma=\bigsqcup_{w\in\Sym_3}\Gamma_w$, where
\begin{align*}
\Gamma_e&=\bigl\{I_3,\ I_3+\vp E_{21},\ I_3+\vp E_{31},\ I_3+\vp E_{32},\ I_3+\vp(E_{21}+E_{32})\bigr\},\\
\Gamma_{(12)}&=\bigl\{P_{(12)},\ (I_3+\vp E_{31})P_{(12)},\ (I_3+\vp E_{32})P_{(12)},\ (I_3+\vp(E_{31}+E_{32}))P_{(12)}\bigr\},\\
\Gamma_{(13)}&=\bigl\{P_{(13)}\bigr\},\\
\Gamma_{(23)}&=\bigl\{P_{(23)},\ (I_3+\vp E_{21})P_{(23)},\ (I_3+\vp E_{31})P_{(23)},\ (I_3+\vp(E_{21}+E_{31}))P_{(23)}\bigr\},\\
\Gamma_{(123)}&=\bigl\{P_{(123)},\ (I_3+\vp E_{32})P_{(123)}\bigr\},\qquad\Gamma_{(132)}=\bigl\{P_{(132)},\ (I_3+\vp E_{21})P_{(132)}\bigr\}.
\end{align*}
\end{proposition}

\begin{proof}
The proposition follows from \cite[Corollary~4.2]{OPV06}.
\end{proof}

\begin{remark}
Note that one can prove the above proposition by using elementary methods of linear algebra and using Lemma~\ref{lem:Bruhat}.
\end{remark}

Recall that $D(\chi)=\dim\End_G(I(\chi))$. For $\delta\in G$ put $H^\delta(I(\chi)):=\Hom_{\delta B\delta^{-1}\cap B}(\chi^\delta,\chi)$, where $\chi^\delta$ is as in~\eqref{eq:conj}. By~\eqref{eq:end},
\begin{equation}\label{eq:Dsum}
D(\chi)=\sum_{\delta\in\Gamma}\dim H^\delta(I(\chi)).
\end{equation}

\begin{proposition}\label{prop:Dchi}
Let $\chi_1,\chi_2,\chi_3\in\Oh$ and let $\{i,j,k\}=\{1,2,3\}$. Then
\[
D(\chi)=\begin{cases}
1 & \text{if } \chi_1,\chi_2,\chi_3 \text{ are pairwise non-equivalent under } \sim,\\
2 & \text{if } \chi_i\neq\chi_j \text{ and } \chi_i\sim\chi_j\not\sim\chi_k,\\
3 & \text{if } \chi_i=\chi_j \text{ and } \chi_i\not\sim\chi_k,\\
7 & \text{if } \chi_1\sim\chi_2\sim\chi_3 \text{ and } \chi_1,\chi_2,\chi_3 \text{ are pairwise distinct},\\
10 & \text{if } \chi_1\sim\chi_2\sim\chi_3 \text{ and } \chi_i=\chi_j\neq\chi_k,\\
18 & \text{if } \chi_1=\chi_2=\chi_3.
\end{cases}
\]
In particular, $I(\chi)$ is irreducible if and only if $\chi_1,\chi_2,\chi_3$ are pairwise non-equivalent under $\sim$.
\end{proposition}

\begin{proof}
Each summand in~\eqref{eq:Dsum} has dimension at most one, so $D(\chi)$ is the number of $\delta\in\Gamma$ with $H^\delta(I(\chi))\neq0$, that is, with $\chi^\delta=\chi$ on $\delta B\delta^{-1}\cap B$. We check this condition for each $\delta\in\Gamma$ by a direct computation. For example, if $\delta=I_3+\vp E_{21}$, then $\delta B\delta^{-1}\cap B$ consists of the $b=(b_{ij})\in B$ with $b_{11}\equiv b_{22}\bmod\vp$, and
\[
\chi^\delta(b)=\chi(b)\,(\chi_1\chi_2^{-1})\bigl(1+\vp b_{11}^{-1}b_{12}\bigr),
\]
so $\chi^\delta=\chi$ exactly when $\chi_1\sim\chi_2$. The other cases are similar, and give the following conditions.
\begin{enumerate}
\item For $\delta=I_3$ there is no condition.
\item For $\delta=I_3+\vp E_{21}$, $I_3+\vp E_{31}$, $I_3+\vp E_{32}$ the condition is $\chi_1\sim\chi_2$, $\chi_1\sim\chi_3$, $\chi_2\sim\chi_3$, respectively, and for $\delta=I_3+\vp(E_{21}+E_{32})$ it is $\chi_1\sim\chi_2\sim\chi_3$.
\item For $\delta=P_{(12)}$, $P_{(13)}$, $P_{(23)}$ the condition is $\chi_1=\chi_2$, $\chi_1=\chi_3$, $\chi_2=\chi_3$, respectively, and for $\delta=P_{(123)}$, $P_{(132)}$ it is $\chi_1=\chi_2=\chi_3$.
\item For $\delta=(I_3+\vp E_{31})P_{(12)}$ and $(I_3+\vp E_{32})P_{(12)}$ the condition is $\chi_1=\chi_2$ and $\chi_1\sim\chi_3$; for $\delta=(I_3+\vp E_{21})P_{(23)}$ and $(I_3+\vp E_{31})P_{(23)}$ it is $\chi_2=\chi_3$ and $\chi_2\sim\chi_1$; and for $\delta=(I_3+\vp E_{32})P_{(123)}$ and $(I_3+\vp E_{21})P_{(132)}$ it is $\chi_1=\chi_3$ and $\chi_1\sim\chi_2$.
\item For $\delta=(I_3+\vp(E_{31}+E_{32}))P_{(12)}$ and $(I_3+\vp(E_{21}+E_{31}))P_{(23)}$ the condition is $\chi_1\sim\chi_2\sim\chi_3$.
\end{enumerate}
Counting, in each case of the proposition, the $\delta$ whose condition holds gives the stated value of $D(\chi)$.
\end{proof}

\section{The cases \texorpdfstring{$D(\chi)\in\{2,3\}$}{D(chi) in \{2,3\}}}\label{sec:D23}

In this section we describe $I(\chi)$ when $D(\chi)\in\{2,3\}$. We use the inclusions $B\subseteq B_1=B\cdot K\subseteq G$ and transitivity of induction,
\[
I(\chi)=\Ind_{B_1}^G\Ind_B^{B_1}\chi.
\]
We first study the irreducible constituents of $\Ind_B^{B_1}\chi$ and then use the value of $D(\chi)$ to determine those of $I(\chi)$.

\begin{lemma}\label{lem:B1dc}
\begin{enumerate}
\item[(1)] $B_1\backslash G/B_1\cong\Sym_3$, and the permutation matrices $\{P_w:w\in\Sym_3\}$ form a set of representatives.
\item[(2)] The set $\Gamma_e=\{I_3,\ I_3+\vp E_{21},\ I_3+\vp E_{31},\ I_3+\vp E_{32},\ I_3+\vp(E_{21}+E_{32})\}$ is a set of representatives for $B\backslash B_1/B$.
\end{enumerate}
\end{lemma}

\begin{proof}
Since $K$ is normal in $G$ and contained in $B_1$, reduction modulo $\vp$ gives a bijection $B_1\backslash G/B_1\to\overline{B}\backslash\GL_3(\Fq)/\overline{B}$, so (1) follows from Lemma~\ref{lem:Bruhat}. Part (2) is the case $w=e$ of Proposition~\ref{prop:18}.
\end{proof}

\begin{lemma}\label{lem:EndB1}
For $\chi_1,\chi_2,\chi_3\in\Oh$,
\[
\dim\End_{B_1}\bigl(\Ind_B^{B_1}\chi\bigr)=\begin{cases}
1&\text{if }\chi_1,\chi_2,\chi_3\text{ are pairwise non-equivalent under }\sim,\\
2&\text{if }\chi_i\sim\chi_j\not\sim\chi_k\text{ for some }\{i,j,k\}=\{1,2,3\},\\
5&\text{if }\chi_1\sim\chi_2\sim\chi_3.
\end{cases}
\]
\end{lemma}

\begin{proof}
This follows from~\eqref{eq:hom}, Lemma~\ref{lem:B1dc}(2) and Proposition~\ref{prop:Dchi}.
\end{proof}

\subsection{The cases \texorpdfstring{$\chi_1\sim\chi_2$}{chi1 ~ chi2} and \texorpdfstring{$\chi_2\sim\chi_3$}{chi2 ~ chi3}}\label{subsec:21-32}
For $1\le j<i\le3$ let $U_{ij}:=\{I_3+\vp xE_{ij}:x\in\cO_2\}$. The following subgroups of $B_1$ will be useful:
\begin{align*}
B_{21}&:=B\cdot U_{21}=\left\{\begin{pmatrix}a&b&c\\ \vp d&e&f\\ 0&0&g\end{pmatrix}:a,e,g\in\cO_2^\times,\ b,c,d,f\in\cO_2\right\},\\
B_{32}&:=B\cdot U_{32}=\left\{\begin{pmatrix}a&b&c\\ 0&d&e\\ 0&\vp f&g\end{pmatrix}:a,d,g\in\cO_2^\times,\ b,c,e,f\in\cO_2\right\}.
\end{align*}
Indeed $B_{21}=\{x\in B_1:x_{31}=x_{32}=0\}$ and $B_{32}=\{x\in B_1:x_{21}=x_{31}=0\}$ are subgroups. We have two chains of subgroups of $B_1$:
\[
B\subseteq B_{21}\subseteq B_1\qquad\text{and}\qquad B\subseteq B_{32}\subseteq B_1.
\]

\begin{lemma}\label{lem:BB21}
The double coset spaces $B\backslash B_{21}/B$ and $B\backslash B_{32}/B$ each consist of exactly two elements, with representatives $\{I_3,\ I_3+\vp E_{21}\}$ and $\{I_3,\ I_3+\vp E_{32}\}$, respectively.
\end{lemma}

\begin{proof}
The proof is elementary and hence is omitted.
\end{proof}

{\chg Let
\[
B_{\delta_1}:=\left\{X=\begin{pmatrix}a&b&c\\ \vp d&e&f\\ 0&0&g\end{pmatrix}\in B_{21}:a-e\in\vp\cO_2\right\}.
\]
This is the kernel of the homomorphism $B_{21}\to\Fq^\times$, $X\mapsto(ae^{-1}\bmod\vp)$; in particular it is a normal subgroup of $B_{21}$ of index $q-1$.}

\begin{proposition}\label{prop:B21}
{\chg Let $\chi_1,\chi_2,\chi_3\in\Oh$ with $\chi_1\sim\chi_2$. For $x=\begin{psmallmatrix}a&b&c\\ \vp d&e&f\\ 0&0&g\end{psmallmatrix}\in B_{21}$ put
\begin{equation}\label{eq:chi21}
\chi^{(21)}(x):=\chi_1(a)\chi_2(e)\chi_3(g)\,\chi_1\bigl(1-\vp a^{-1}e^{-1}bd\bigr)=\chi_1(ae-\vp bd)\,(\chi_1^{-1}\chi_2)(e)\,\chi_3(g),
\end{equation}
and for $X\in B_{\delta_1}$, written in the same way, put
\[
\chi^{\delta_1,\psi}(X):=\chi^{(21)}(X)\,\psi(\vp e^{-1}d).
\]
Then $\chi^{(21)}$ is a one-dimensional character of $B_{21}$ extending $\chi$, $\chi^{\delta_1,\psi}$ is a character of $B_{\delta_1}$ such that $\chi^{\delta_1,\psi}=\chi$ on $B_{\delta_1}\cap B$, and
\[
\Ind_B^{B_{21}}\chi\;\cong\;\chi^{(21)}\oplus\tau^{(21)},\qquad\text{where }\ \tau^{(21)}:=\Ind_{B_{\delta_1}}^{B_{21}}\chi^{\delta_1,\psi}
\]
is an irreducible representation of $B_{21}$ of dimension $q-1$.}
\end{proposition}

\begin{proof}
Note that $ae-\vp bd$ is the determinant of the upper left $2\times2$ block of $x$, which is multiplicative on $B_{21}$, and $x\mapsto e\bmod\vp$ is multiplicative on $B_{21}\subseteq B_1$. Since $\chi_1^{-1}\chi_2$ is trivial on $1+\vp\cO_2$, the value $(\chi_1^{-1}\chi_2)(e)$ depends only on $e\bmod\vp$. Hence $\chi^{(21)}$ is a character of $B_{21}$, and clearly $\chi^{(21)}|_B=\chi$.

{\chg For $X,X'\in B_{\delta_1}$ we have $(XX')_{21}=\vp(da'+ed')$ and $(XX')_{22}\equiv ee'\bmod\vp$; since $a'\equiv e'\bmod\vp$,
\[
\frac{da'+ed'}{ee'}=\frac{d}{e}\cdot\frac{a'}{e'}+\frac{d'}{e'}\equiv\frac{d}{e}+\frac{d'}{e'}\pmod{\vp},
\]
so $\theta\colon X\mapsto\psi(\vp e^{-1}d)$ is a character of $B_{\delta_1}$. Hence $\chi^{\delta_1,\psi}=\chi^{(21)}\theta$ is a character, and it agrees with $\chi$ on $B_{\delta_1}\cap B$, where $d=0$.

We now apply Lemma~\ref{lem:split} with $H=B$, $H^+=B_{21}$, $\lambda^+=\chi^{(21)}$, $H'=B_{\delta_1}$ and $\theta$, so that $\lambda^+|_{H'}\theta=\chi^{\delta_1,\psi}$. Since $d$ only matters modulo $\vp$, we have $[B_{21}:B]=q$. The elements $\gamma(y):=\diag(y,1,1)\in B$, $y\in\cO_2^\times$, represent the cosets of $B_{\delta_1}$ in $B_{21}$; hence $B_{21}=B\,B_{\delta_1}$. The character $\theta$ is trivial on $B\cap B_{\delta_1}$, and it is non-trivial because $\psi|_{\vp\cO_2}\neq\one$. Finally, $\theta^{\gamma(y)}(X)=\theta(\gamma(y)^{-1}X\gamma(y))=\psi(\vp e^{-1}dy)$, so $\theta^{\gamma(y)}=\theta$ if and only if $y\in1+\vp\cO_2$, that is, if and only if $\gamma(y)\in B_{\delta_1}$. As $\theta$ is a character of the normal subgroup $B_{\delta_1}$, the character $\theta^x$ only depends on the coset $xB_{\delta_1}$, and therefore $\theta^x\neq\theta$ for every $x\in B_{21}\smallsetminus B_{\delta_1}$. Lemma~\ref{lem:split} now gives the proposition.}
\end{proof}

\begin{corollary}\label{cor:21}
Let $\chi_1,\chi_2,\chi_3\in\Oh$ be pairwise distinct with $\chi_1\sim\chi_2\not\sim\chi_3$. Then
\[
I(\chi)\;\cong\;\Ind_{B_{21}}^G\chi^{(21)}\oplus\Ind_{B_{\delta_1}}^G\chi^{\delta_1,\psi},
\]
and the two summands are irreducible and non-isomorphic, of dimensions $q^2(q+1)(q^2+q+1)$ and $q^2(q^2-1)(q^2+q+1)$.
\end{corollary}

\begin{proof}
The corollary follows from transitivity of induction and using Proposition~\ref{prop:B21}.
\end{proof}

\begin{proposition}\label{prop:B32}
{\chg Let $\chi_1,\chi_2,\chi_3\in\Oh$ with $\chi_2\sim\chi_3$. For $x=\begin{psmallmatrix}a&b&c\\ 0&d&e\\ 0&\vp f&g\end{psmallmatrix}\in B_{32}$ put
\begin{equation}\label{eq:chi32}
\chi^{(32)}(x):=\chi_1(a)\chi_2(d)\chi_3(g)\,\chi_2\bigl(1-\vp d^{-1}g^{-1}ef\bigr)=\chi_1(a)\,\chi_2(dg-\vp ef)\,(\chi_2^{-1}\chi_3)(g).
\end{equation}
Let
\[
B_{32}^\circ:=\{x\in B_{32}:d-g\in\vp\cO_2\},
\]
and for $x\in B_{32}^\circ$ put $\chi^{(32),\psi}(x):=\chi^{(32)}(x)\,\psi(\vp g^{-1}f)$. 
Then $\chi^{(32)}$ is a one-dimensional character of $B_{32}$ extending $\chi$, $\chi^{(32),\psi}$ is a character of $B_{32}^\circ$, and
\[
\Ind_B^{B_{32}}\chi\;\cong\;\chi^{(32)}\oplus\tau^{(32)},\qquad\text{where }\ \tau^{(32)}:=\Ind_{B_{32}^\circ}^{B_{32}}\chi^{(32),\psi}
\]
is an irreducible representation of $B_{32}$ of dimension $q-1$.}
\end{proposition}

\begin{proof}
{\chg As in the proof of Proposition~\ref{prop:B21}, using that $dg-\vp ef$ is the determinant of the lower right $2\times2$ block of $x$, one sees that $\chi^{(32)}$ is a character of $B_{32}$ extending $\chi$. The subgroup $B_{32}^\circ$ is the kernel of the homomorphism $B_{32}\to\Fq^\times$, $x\mapsto dg^{-1}\bmod\vp$, so it is normal of index $q-1$, and its cosets are represented by the elements $\diag(1,1,y)\in B$, $y\in\cO_2^\times$. For $x,x'\in B_{32}^\circ$ we have $(xx')_{32}=\vp(fd'+gf')$ and $(xx')_{33}\equiv gg'\bmod\vp$, and $d'\equiv g'\bmod\vp$; hence $\theta\colon x\mapsto\psi(\vp g^{-1}f)$ is a character of $B_{32}^\circ$, and so is $\chi^{(32),\psi}=\chi^{(32)}\theta$. Now apply Lemma~\ref{lem:split} with $H=B$, $H^+=B_{32}$, $\lambda^+=\chi^{(32)}$, $H'=B_{32}^\circ$ and $\theta$: the character $\theta$ is trivial on $B\cap B_{32}^\circ$, where $f=0$, it is non-trivial, and $\diag(1,1,y)$ changes it to $x\mapsto\psi(\vp g^{-1}fy^{-1})$. The remaining hypotheses are checked as in the proof of Proposition~\ref{prop:B21}.}
\end{proof}

\begin{corollary}\label{cor:32}
Let $\chi_1,\chi_2,\chi_3\in\Oh$ be pairwise distinct with $\chi_2\sim\chi_3\not\sim\chi_1$. Then
\[
I(\chi)\;\cong\;\Ind_{B_{32}}^G\chi^{(32)}\oplus\Ind_{B_{32}^\circ}^{G}\chi^{(32),\psi},
\]
and the two summands are irreducible and non-isomorphic.
\end{corollary}

\begin{proof}
The corollary follows from transitivity of induction and Proposition~\ref{prop:B32}.
\end{proof}

\subsection{The case \texorpdfstring{$\chi_1\sim\chi_3\not\sim\chi_2$}{chi1 ~ chi3}}\label{subsec:13}
The subset $B\cdot U_{31}$ of $B_1$ is not a subgroup. Consequently, if $\chi_1\sim\chi_3$, the argument of Propositions~\ref{prop:B21} and~\ref{prop:B32} does not apply directly. To overcome this, we consider the following character $\chi_{(12)}:=\chi_2\otimes\chi_1\otimes\chi_3$.

\begin{lemma}\label{lem:swap}
Let $\chi_1,\chi_2,\chi_3\in\Oh$.
\begin{enumerate}
\item[(1)] If $\chi_1\sim\chi_3$, then $\dim\Hom_G\bigl(I(\chi),I(\chi_{(12)})\bigr)\ge2$.
\item[(2)] If $\chi_1,\chi_2,\chi_3$ are pairwise distinct and $\chi_1\sim\chi_3\not\sim\chi_2$, then $I(\chi)\cong I(\chi_{(12)})$.
\end{enumerate}
\end{lemma}

\begin{proof}
(1) By~\eqref{eq:hom},
\[
\Hom_G\bigl(I(\chi),I(\chi_{(12)})\bigr)\cong\bigoplus_{s\in[B\backslash G/B]}\Hom_{sBs^{-1}\cap B}\bigl(\chi^s,\chi_{(12)}\bigr).
\]
It is easy to check that for $s=P_{(12)}$ and $s=(I_3+\vp E_{32})P_{(12)}$, $\Hom_{sBs^{-1}\cap B}\bigl(\chi^s,\chi_{(12)}\bigr)\neq0$ as for $s=P_{(12)}$ we have $\chi^s=\chi_{(12)}$ on $sBs^{-1}\cap B$ for every $\chi$, and for $s=(I_3+\vp E_{32})P_{(12)}$ this holds because $\chi_1\sim\chi_3$. Hence the dimension is at least two.

(2) By Proposition~\ref{prop:Dchi}, $\dim\End_G(I(\chi))=2$. In fact, Proposition~\ref{prop:Dchi} gives $\dim\End_G(I(\chi_{(12)}))=2$. Hence $I(\chi)$ and $I(\chi_{(12)})$ are each a direct sum of two non-isomorphic irreducible representations. Therefore $\dim\Hom_G\bigl(I(\chi),I(\chi_{(12)})\bigr)$ is the number of irreducible constituents they have in common, which is at most two. By (1) it equals two, so $I(\chi)$ and $I(\chi_{(12)})$ have the same constituents, and $I(\chi)\cong I(\chi_{(12)})$.
\end{proof}

Let $B^*:=P_{(12)}BP_{(12)}$ and let $P_{1,2}$ be the subgroup of $G$ of block upper triangular matrices with diagonal blocks of sizes $1\times1$ and $2\times2$. Put $P^*:=P_{(12)}P_{1,2}P_{(12)}$. Explicitly,
\[
B^*=\left\{\begin{pmatrix}a&0&b\\ c&d&e\\ 0&0&f\end{pmatrix}:\begin{array}{l}a,d,f\in\cO_2^\times,\\ b,c,e\in\cO_2\end{array}\right\},\qquad
P^*=\left\{\begin{pmatrix}a&0&b\\ c&d&e\\ f&0&g\end{pmatrix}:\begin{array}{l}ag-bf,\ d\in\cO_2^\times,\\ c,e\in\cO_2\end{array}\right\}.
\]
Let $P^*_\vp:=\{X\in P^*:X_{31}\in\vp\cO_2\}$. Then $B^*\subseteq P^*_\vp$ and $[P^*_\vp:B^*]=q$; in fact $P^*_\vp=P_{(12)}B_{32}P_{(12)}$. Define a character $\chi^*$ of $B^*$ by
\[
\chi^*\begin{pmatrix}a&0&b\\ c&d&e\\ 0&0&f\end{pmatrix}:=\chi_{(12)}\begin{pmatrix}d&c&e\\ 0&a&b\\ 0&0&f\end{pmatrix}=\chi_1(a)\chi_2(d)\chi_3(f),
\]
that is, $\chi^*(X)=\chi_{(12)}(P_{(12)}XP_{(12)})$. Then $\Ind_{B^*}^G\chi^*\cong I(\chi_{(12)})$, and hence, by Lemma~\ref{lem:swap}, $I(\chi)\cong\Ind_{B^*}^G\chi^*$ whenever $\chi_1,\chi_2,\chi_3$ are pairwise distinct and $\chi_1\sim\chi_3\not\sim\chi_2$. It therefore suffices to describe $\Ind_{B^*}^G\chi^*$ in this case.

\begin{proposition}\label{prop:Pstar}
{\chg Let $\chi_1,\chi_2,\chi_3\in\Oh$ with $\chi_1\sim\chi_3$. For $X=\begin{psmallmatrix}a&0&b\\ c&d&e\\ \vp f&0&g\end{psmallmatrix}\in P^*_\vp$ put
\[
\chi^\vp(X):=\chi_1(a)\chi_2(d)\chi_3(g)\,\chi_1\bigl(1-\vp a^{-1}g^{-1}bf\bigr)=\chi_1(ag-\vp bf)\,\chi_2(d)\,(\chi_1^{-1}\chi_3)(g),
\]
let $P^{*\circ}_\vp:=\{X\in P^*_\vp:a-g\in\vp\cO_2\}$, and for $X\in P^{*\circ}_\vp$ put $\chi^{\vp,\psi}(X):=\chi^\vp(X)\,\psi(\vp g^{-1}f)$. Then $\chi^\vp$ is a one-dimensional character of $P^*_\vp$ extending $\chi^*$, $\chi^{\vp,\psi}$ is a character of $P^{*\circ}_\vp$, and
\[
\Ind_{B^*}^{P^*_\vp}\chi^*\;\cong\;\chi^\vp\oplus\tau^\vp,\qquad\text{where }\ \tau^\vp:=\Ind_{P^{*\circ}_\vp}^{P^*_\vp}\chi^{\vp,\psi}
\]
is an irreducible representation of $P^*_\vp$ of dimension $q-1$.}
\end{proposition}

\begin{proof}
{\chg Conjugation $X\mapsto P_{(12)}XP_{(12)}$ maps $B^*$ onto $B$, $P^*_\vp$ onto $B_{32}$ and $P^{*\circ}_\vp$ onto $B_{32}^\circ$, and it carries $\chi^*$ to $\chi_{(12)}$, $\chi^\vp$ to the character $(\chi_{(12)})^{(32)}$ of~\eqref{eq:chi32}, and $\chi^{\vp,\psi}$ to $(\chi_{(12)})^{(32),\psi}$. The second and third components $\chi_1,\chi_3$ of $\chi_{(12)}$ are equivalent, so Proposition~\ref{prop:B32} applies to $\chi_{(12)}$, and transporting it back by $P_{(12)}$ gives the claim.}
\end{proof}

\begin{corollary}\label{cor:13}
Let $\chi_1,\chi_2,\chi_3\in\Oh$ be pairwise distinct with $\chi_1\sim\chi_3\not\sim\chi_2$. Then
\[
I(\chi)\;\cong\;\Ind_{P^*_\vp}^G\chi^\vp\oplus\Ind_{P^{*\circ}_\vp}^G\chi^{\vp,\psi},
\]
and the two summands are irreducible and non-isomorphic.
\end{corollary}

\begin{proof}
By Lemma~\ref{lem:swap}, $I(\chi)\cong\Ind_{B^*}^G\chi^*=\Ind_{P^*_\vp}^G\Ind_{B^*}^{P^*_\vp}\chi^*$. The claim follows from Proposition~\ref{prop:Pstar} and $D(\chi)=2$.
\end{proof}

\subsection{The case \texorpdfstring{$D(\chi)=3$}{D(chi)=3}}
Let $P_{2,1}$ be the subgroup of $G$ of block upper triangular matrices with diagonal blocks of sizes $2\times2$ and $1\times1$. For a representation $\pi$ of $\GL_2(\cO_2)$ and $\lambda\in\Oh$ we write $\pi\boxtimes\lambda$ for the representation $\begin{psmallmatrix}p_{11}&p_{12}\\ 0&p_{22}\end{psmallmatrix}\mapsto\pi(p_{11})\lambda(p_{22})$ of $P_{2,1}$, where $p_{11}\in\GL_2(\cO_2)$ and $p_{22}\in\cO_2^\times$. Let $\St_2$ be the Steinberg representation of $\GL_2(\Fq)$, inflated to $\GL_2(\cO_2)$.

\begin{theorem}\label{thm:D3}
Let $\chi_1,\chi_2,\chi_3\in\Oh$ with $\chi_1=\chi_2\not\sim\chi_3$. Then
\[
I(\chi)\;\cong\;\Ind_{P_{2,1}}^G\bigl((\chi_1\circ\det)\boxtimes\chi_3\bigr)\;\oplus\;\Ind_{P_{2,1}}^G\bigl(((\chi_1\circ\det)\otimes\St_2)\boxtimes\chi_3\bigr)\;\oplus\;\Ind_{B_{21}}^G\tau^{(21)},
\]
and the three summands are irreducible and pairwise non-isomorphic, of dimensions $q^2(q^2+q+1)$, $q^3(q^2+q+1)$ and $q^2(q^2-1)(q^2+q+1)$.
\end{theorem}

\begin{proof}
By Proposition~\ref{prop:B21},
\begin{equation}\label{eq:D3a}
\Ind_B^{B_{21}}\chi\cong\chi^{(21)}\oplus\tau^{(21)}.
\end{equation}
For $x=\begin{psmallmatrix}A&*\\ 0&g\end{psmallmatrix}\in B_{21}$ with $A=\begin{psmallmatrix}a&b\\ \vp d&e\end{psmallmatrix}$ we have $\det A=ae-\vp bd$. Since $\chi_1=\chi_2$, \eqref{eq:chi21} gives $\chi^{(21)}(x)=\chi_1(\det A)\chi_3(g)$, i.e., $\chi^{(21)}=\eta|_{B_{21}}$ and $\eta:=(\chi_1\circ\det)\boxtimes\chi_3$. Hence
\begin{equation}\label{eq:D3b}
\Ind_{B_{21}}^{P_{2,1}}\chi^{(21)}=\Ind_{B_{21}}^{P_{2,1}}\bigl(\eta|_{B_{21}}\bigr)\cong\eta\otimes\Ind_{B_{21}}^{P_{2,1}}\one.
\end{equation}
Since $P_{2,1}\cap K\subseteq B_{21}$, the representation $\Ind_{B_{21}}^{P_{2,1}}\one$ factors through $P_{2,1}/(P_{2,1}\cap K)\cong\overline{P}_{2,1}$, and hence $\Ind_{B_{21}}^{P_{2,1}}\one$ should be realised as inflation of the representation $\Ind_{\overline{B}}^{\overline{P}_{2,1}}\one$ to $P_{2,1}$. Let $\overline{B}_2$ be the upper triangular subgroup of $\GL_2(\Fq)$. Since $\overline{P}_{2,1}=(\GL_2(\Fq)\times\Fq^\times)\ltimes N$, where $N$ is the unipotent radical of $\overline{P}_{2,1}$, and $\overline{B}=(\overline{B}_2\times\Fq^\times)\ltimes N$, we have $\Ind_{\overline{B}}^{\overline{P}_{2,1}}\one\cong\bigl(\Ind_{\overline{B}_2}^{\GL_2(\Fq)}\one\bigr)\boxtimes\one$. Using $\Ind_{\overline{B}_2}^{\GL_2(\Fq)}\one\cong\one\oplus\St_2$ and~\eqref{eq:D3b}, we obtain
\[
\Ind_{B_{21}}^{P_{2,1}}\chi^{(21)}\;\cong\;\bigl((\chi_1\circ\det)\boxtimes\chi_3\bigr)\oplus\bigl(((\chi_1\circ\det)\otimes\St_2)\boxtimes\chi_3\bigr),
\]
where $\St$ is inflated $\GL_2(\cO_2)$ via the natural map $\GL_2(\cO_2)\twoheadrightarrow\GL_2(\Fq)$. By transitivity of induction and~\eqref{eq:D3a}, $I(\chi)$ is the stated direct sum of three non-zero representations. Since $D(\chi)=3$ the theorem follows from Proposition~\ref{prop:Dchi}.
\end{proof}

\begin{remark}\label{rem:D3}
Analogous statements hold in the two remaining cases with $D(\chi)=3$. If $\chi_2=\chi_3\not\sim\chi_1$, one uses $P_{1,2}$, $B_{32}$ and $\tau^{(32)}$ instead of $P_{2,1}$, $B_{21}$ and $\tau^{(21)}$. If $\chi_1=\chi_3\not\sim\chi_2$, one uses $P^*$, $B^*$ and $\tau^\vp$ instead of $P_{2,1}$, $B_{21}$ and $\tau^{(21)}$. We omit the details.
\end{remark}

\section{The case \texorpdfstring{$D(\chi)=7$}{D(chi)=7}}\label{sec:D7}

In this section and the next we assume that $\chi_1\sim\chi_2\sim\chi_3$. The main result of this section, Theorem~\ref{thm:D7}, describes $I(\chi)$ when moreover $\chi_1,\chi_2,\chi_3$ are pairwise distinct, so that $D(\chi)=7$. Most of the intermediate results hold for all $\chi$ with $\chi_1\sim\chi_2\sim\chi_3$ and are used again in Section~\ref{sec:D1018}.

Recall the subgroups $B_{21}$ and $B_{32}$ from Section~\ref{subsec:21-32}, and put
\[
B_{21,32}:=B_{21}B_{32}=\{x\in B_1:x_{31}=0\}=\left\{\begin{pmatrix}a&b&c\\ \vp d&e&f\\ 0&\vp h&g\end{pmatrix}:a,e,g\in\cO_2^\times,\ b,c,d,f,h\in\cO_2\right\}.
\]
The second description shows that $B_{21,32}$ is a subgroup. We obtain the chains of subgroups
\[
B\subset B_{21}\subset B_{21,32}\subset B_1\subset G,\qquad B\subset B_{32}\subset B_{21,32}\subset B_1\subset G.
\]

\begin{lemma}\label{lem:chitilde}
Let $\chi_1,\chi_2,\chi_3\in\Oh$ with $\chi_1\sim\chi_2\sim\chi_3$. Define $\ct\colon B_1=B\cdot K\to\C^\times$ by
\[
\ct(bk):=\chi(b)\,(\chi_1\circ\det)(k)\qquad(b\in B,\ k\in K).
\]
Then:
\begin{enumerate}
\item[(1)] $\ct$ is a well-defined character of $B_1$; explicitly, $\ct(x)=(\chi_1\circ\det)(x)\,(\chi_1^{-1}\chi_2)(x_{22})\,(\chi_1^{-1}\chi_3)(x_{33})$ for $x\in B_1$;
\item[(2)] $\ct|_{B_{21}}=\chi^{(21)}$ and $\ct|_{B_{32}}=\chi^{(32)}$.
\end{enumerate}
\end{lemma}

\begin{proof}
The map $x\mapsto x_{ii}\bmod\vp$ is a homomorphism $B_1\to\Fq^\times$, and $\chi_1^{-1}\chi_2$ and $\chi_1^{-1}\chi_3$ are trivial on $1+\vp\cO_2$. Hence the right hand side of the formula in (1) is a character of $B_1$. For $b\in B$ it equals $\chi(b)$, and for $k\in K$ it equals $\chi_1(\det k)$, so it agrees with $\ct$; in particular $\ct$ is well defined. Part (2) follows by comparing with~\eqref{eq:chi21} and~\eqref{eq:chi32}.
\end{proof}

The following lemma describes several double cosets that will be used later; the verification is straightforward and we omit it.

\begin{lemma}\label{lem:dc5}
\begin{enumerate}
\item[(1)] The set $\{I_3,\ I_3+\vp E_{31},\ I_3+\vp E_{32}\}$ is a set of representatives for $B_{21}\backslash B_1/B_{21}$, and also for $B_{21}\backslash B_1/B$.
\item[(2)] The set $\{I_3,\ I_3+\vp E_{21},\ I_3+\vp E_{31}\}$ is a set of representatives for $B_{32}\backslash B_1/B_{32}$, and also for $B_{32}\backslash B_1/B$.
\end{enumerate}
\end{lemma}

\begin{lemma}\label{lem:End3}
Let $\chi_1,\chi_2,\chi_3\in\Oh$ with $\chi_1\sim\chi_2\sim\chi_3$. Then
\[
\dim\End_{B_1}\bigl(\Ind_{B_{21}}^{B_1}\chi^{(21)}\bigr)=\dim\Hom_{B_1}\bigl(\Ind_B^{B_1}\chi,\ \Ind_{B_{21}}^{B_1}\chi^{(21)}\bigr)=3.
\]
\end{lemma}

\begin{proof}
The proof is straightforward and hence left.
\end{proof}

Define $\chi^{(21,32)}:=\ct|_{B_{21,32}}$. By Lemma~\ref{lem:chitilde}, $\chi^{(21,32)}|_{B_{21}}=\chi^{(21)}$ and $\chi^{(21,32)}|_{B_{32}}=\chi^{(32)}$. {\chg Let
\[
B_{\delta_2}:=\left\{Y=\begin{pmatrix}a&b&c\\ \vp d&e&f\\ 0&\vp h&g\end{pmatrix}\in B_{21,32}:e-g\in\vp\cO_2\right\},
\]
a normal subgroup of $B_{21,32}$ of index $q-1$, and if $\chi_1\sim\chi_2\sim\chi_3$ define $\chi^{\delta_2,\psi}\colon B_{\delta_2}\to\C^\times$ by $\chi^{\delta_2,\psi}(Y):=\chi^{(21,32)}(Y)\,\psi(\vp g^{-1}h)$.}

\begin{lemma}\label{lem:B2132}
Let $\chi_1,\chi_2,\chi_3\in\Oh$ with $\chi_1\sim\chi_2$. If $\chi_1\sim\chi_2\sim\chi_3$, then $\chi^{\delta_2,\psi}$ is a character of $B_{\delta_2}$ such that $\chi^{\delta_2,\psi}=\chi^{(21)}$ on $B_{\delta_2}\cap B_{21}$, and
\[
\Ind_{B_{21}}^{B_{21,32}}\chi^{(21)}\;\cong\;\chi^{(21,32)}\oplus\tau^{(21,32)},\qquad\text{where }\ \tau^{(21,32)}:=\Ind_{B_{\delta_2}}^{B_{21,32}}\chi^{\delta_2,\psi}
\]
is an irreducible representation of $B_{21,32}$ of dimension $q-1$.
\end{lemma}

\begin{proof}
{\chg The second case is elementary and we omit it. Suppose that $\chi_1\sim\chi_2\sim\chi_3$. For $Y,Y'\in B_{\delta_2}$ we have $(YY')_{32}=\vp(he'+gh')$ and $(YY')_{33}\equiv gg'\bmod\vp$; since $e'\equiv g'\bmod\vp$, the map $\theta\colon Y\mapsto\psi(\vp g^{-1}h)$ is a character of $B_{\delta_2}$, as in the proof of Proposition~\ref{prop:B21}. Hence $\chi^{\delta_2,\psi}=\chi^{(21,32)}\theta$ is a character, and it agrees with $\chi^{(21)}$ on $B_{\delta_2}\cap B_{21}$, where $h=0$.

We apply Lemma~\ref{lem:split} with $H=B_{21}$, $H^+=B_{21,32}$, $\lambda^+=\chi^{(21,32)}$, $H'=B_{\delta_2}$ and $\theta$, so that $\lambda^+|_H=\chi^{(21)}$ and $\lambda^+|_{H'}\theta=\chi^{\delta_2,\psi}$. Here $[B_{21,32}:B_{21}]=q$, since $h$ only matters modulo $\vp$. The subgroup $B_{\delta_2}$ is the kernel of the surjective homomorphism $B_{21,32}\to\Fq^\times$, $Y\mapsto eg^{-1}\bmod\vp$, and the elements $\gamma(x):=\diag(1,x,1)\in B_{21}$, $x\in\cO_2^\times$, represent its cosets; hence $B_{21,32}=B_{21}B_{\delta_2}$. The character $\theta$ is trivial on $B_{21}\cap B_{\delta_2}$ and non-trivial. Finally $\theta^{\gamma(x)}(Y)=\psi(\vp g^{-1}hx)$, so $\theta^{\gamma(x)}=\theta$ if and only if $x\in1+\vp\cO_2$; as in the proof of Proposition~\ref{prop:B21} it follows that $\theta^x\neq\theta$ for every $x\in B_{21,32}\smallsetminus B_{\delta_2}$. Lemma~\ref{lem:split} now gives the claim.}
\end{proof}

{\chg Let
\[
B_{\delta_3}:=\left\{Z=\begin{pmatrix}a&b&c\\ \vp d&e&f\\ \vp k&\vp h&g\end{pmatrix}\in B_1:a-g\in\vp\cO_2\right\},
\]
a normal subgroup of $B_1$ of index $q-1$, and if $\chi_1\sim\chi_2\sim\chi_3$ define $\chi^{\delta_3,\psi}\colon B_{\delta_3}\to\C^\times$ by $\chi^{\delta_3,\psi}(Z):=\ct(Z)\,\psi(\vp g^{-1}k)$.}

\begin{lemma}\label{lem:tildetau}
{\chg Let $\chi_1,\chi_2,\chi_3\in\Oh$ with $\chi_1\sim\chi_2\sim\chi_3$. Then $\chi^{\delta_3,\psi}$ is a character of $B_{\delta_3}$ such that $\chi^{\delta_3,\psi}=\chi^{(21,32)}$ on $B_{\delta_3}\cap B_{21,32}$, and
\[
\Ind_{B_{21,32}}^{B_1}\chi^{(21,32)}\;\cong\;\ct\oplus\widetilde{\tau},\qquad\text{where }\ \widetilde{\tau}:=\Ind_{B_{\delta_3}}^{B_1}\chi^{\delta_3,\psi}
\]
is an irreducible representation of $B_1$ of dimension $q-1$.}
\end{lemma}

\begin{proof}
{\chg For $Z,Z'\in B_{\delta_3}$ we have $(ZZ')_{31}=\vp(ka'+gk')$ and $(ZZ')_{33}\equiv gg'\bmod\vp$. Since $a'\equiv g'\bmod\vp$, the map $\theta\colon Z\mapsto\psi(\vp g^{-1}k)$ is a character of $B_{\delta_3}$. Hence $\chi^{\delta_3,\psi}=\ct\,\theta$ is a character, and it agrees with $\chi^{(21,32)}$ on $B_{\delta_3}\cap B_{21,32}$, where $k=0$.

We apply Lemma~\ref{lem:split} with $H=B_{21,32}$, $H^+=B_1$, $\lambda^+=\ct$, $H'=B_{\delta_3}$ and $\theta$, so that $\lambda^+|_H=\chi^{(21,32)}$ and $\lambda^+|_{H'}\theta=\chi^{\delta_3,\psi}$. Here $[B_1:B_{21,32}]=q$, since $k$ only matters modulo $\vp$. The subgroup $B_{\delta_3}$ is the kernel of the surjective homomorphism $B_1\to\Fq^\times$, $Z\mapsto ag^{-1}\bmod\vp$, and the elements $\gamma(x):=\diag(1,1,x)\in B_{21,32}$, $x\in\cO_2^\times$, represent its cosets; hence $B_1=B_{21,32}B_{\delta_3}$. The character $\theta$ is trivial on $B_{21,32}\cap B_{\delta_3}$ and non-trivial. Finally $\theta^{\gamma(x)}(Z)=\psi(\vp g^{-1}kx^{-1})$, so $\theta^{\gamma(x)}=\theta$ if and only if $x\in1+\vp\cO_2$, and as before $\theta^x\neq\theta$ for every $x\in B_1\smallsetminus B_{\delta_3}$. Lemma~\ref{lem:split} now gives the claim.}
\end{proof}

\begin{lemma}\label{lem:B1sum}
Let $\chi_1,\chi_2,\chi_3\in\Oh$ with $\chi_1\sim\chi_2\sim\chi_3$.
Then
\begin{enumerate}
\item 
\(\Ind_{B_{21}}^{B_1}\chi^{(21)}\cong\ct\oplus\widetilde{\tau}\oplus\Ind_{B_{21,32}}^{B_1}\tau^{(21,32)}\).
\item
\(\Ind_B^{B_1}\chi\cong\ct\oplus\widetilde{\tau}\oplus\Ind_{B_{21,32}}^{B_1}\tau^{(21,32)}\oplus\Ind_{B_{21}}^{B_1}\tau^{(21)}\).
\end{enumerate}
\end{lemma}

\begin{proof}
By transitivity of induction and Lemma~\ref{lem:B2132} and Lemma~\ref{lem:tildetau},
\[
\Ind_{B_{21}}^{B_1}\chi^{(21)}\cong\Ind_{B_{21,32}}^{B_1}\chi^{(21,32)}\oplus\Ind_{B_{21,32}}^{B_1}\tau^{(21,32)}\cong\ct\oplus\widetilde{\tau}\oplus\Ind_{B_{21,32}}^{B_1}\tau^{(21,32)}.
\]
The second isomorphism follows from the first and $\Ind_B^{B_1}\chi\cong\Ind_{B_{21}}^{B_1}\chi^{(21)}\oplus\Ind_{B_{21}}^{B_1}\tau^{(21)}$, see Proposition~\ref{prop:B21}.
\qedhere 
\end{proof}

\begin{lemma}\label{lem:B1irr}
Let $\chi_1,\chi_2,\chi_3\in\Oh$ with $\chi_1\sim\chi_2\sim\chi_3$.
Then:
\begin{enumerate}
\item[(1)] the representations $\ct$, $\widetilde{\tau}$ and $\Ind_{B_{21,32}}^{B_1}\tau^{(21,32)}$ of $B_1$ are irreducible and pairwise non-isomorphic;
\item[(2)] the representation $\Ind_{B_{21}}^{B_1}\tau^{(21)}$ is the direct sum of two non-isomorphic irreducible representations, neither of which is isomorphic to a representation in (1).
\end{enumerate}
\end{lemma}

\begin{proof}
(1) By Lemma~\ref{lem:B1sum}, $V:=\Ind_{B_{21}}^{B_1}\chi^{(21)}$ is the direct sum of the three representations in (1), and $\dim\End_{B_1}(V)=3$ by Lemma~\ref{lem:End3}.

(2) Let $W:=\Ind_{B_{21}}^{B_1}\tau^{(21)}$, so that $\Ind_B^{B_1}\chi\cong V\oplus W$. Now, (2) follows from $\dim\End_{B_1}(V\oplus W)=5$ by Lemma~\ref{lem:EndB1}.
\end{proof}

{\chg Combining Lemmas~\ref{lem:B1sum} and~\ref{lem:B1irr} with the descriptions of $\tau^{(21)}$, $\tau^{(21,32)}$ and $\widetilde{\tau}$ in Proposition~\ref{prop:B21} and Lemmas~\ref{lem:B2132} and~\ref{lem:tildetau}, we obtain the following theorem.}

\begin{theorem}\label{thm:B1}
Let $\chi_1,\chi_2,\chi_3\in\Oh$ with $\chi_1\sim\chi_2\sim\chi_3$, and let $B_{\delta_1},B_{\delta_2},B_{\delta_3}$ be as above. Then
\[
\Ind_B^{B_1}\chi\;\cong\;\ct\;\oplus\;\Ind_{B_{\delta_3}}^{B_1}\chi^{\delta_3,\psi}\;\oplus\;\Ind_{B_{\delta_2}}^{B_1}\chi^{\delta_2,\psi}\;\oplus\;\Ind_{B_{\delta_1}}^{B_1}\chi^{\delta_1,\psi},
\]
where $\ct$, $\Ind_{B_{\delta_3}}^{B_1}\chi^{\delta_3,\psi}$ and $\Ind_{B_{\delta_2}}^{B_1}\chi^{\delta_2,\psi}$ are irreducible, and $\Ind_{B_{\delta_1}}^{B_1}\chi^{\delta_1,\psi}$ is the direct sum of two non-isomorphic irreducible representations.
\end{theorem}

The dimensions of the four summands are $1$, $q-1$, $q(q-1)$ and $q^2(q-1)$, which add up to $q^3=[B_1:B]$. By Lemma~\ref{lem:B1irr}, the representation $\Ind_{B_{\delta_1}}^{B_1}\chi^{\delta_1,\psi}$ has exactly two irreducible components. We describe them in the following subsection.

\subsection{The two components of \texorpdfstring{$\Ind_{B_{\delta_1}}^{B_1}\chi^{\delta_1,\psi}$}{Ind(B delta1,B1)}}
Define the following two subgroups of $B_1$:
\begin{align*}
B_{\delta_1}^{(1)}&:=\left\{X=\begin{pmatrix}a&b&c\\ \vp d&e&f\\ 0&\vp h&g\end{pmatrix}\in B_{21,32}:a-e\in\vp\cO_2\right\},\\
B_{\delta_1}^{(2)}&:=\bigl\{X\in B_{\delta_1}^{(1)}:a-g\in\vp\cO_2\bigr\}.
\end{align*}

\begin{proposition}\label{prop:twocomp}
Let $\chi_1,\chi_2,\chi_3\in\Oh$ with $\chi_1\sim\chi_2\sim\chi_3$, and let $\psi$ be a character of $\cO_2$ with $\psi|_{\vp\cO_2}\neq\one$. Define $\chi^{(1,\psi)}_{\delta_1}\colon B_{\delta_1}^{(1)}\to\C^\times$ by
\[
\chi^{(1,\psi)}_{\delta_1}(X):=\ct(X)\,\psi(\vp e^{-1}d)=\chi_1(a)\chi_2(e)\chi_3(g)\,\chi_1\bigl(1-\vp a^{-1}e^{-1}bd-\vp e^{-1}g^{-1}fh\bigr)\,\psi(\vp e^{-1}d),
\]
and $\chi^{(2,\psi)}_{\delta_1}\colon B_{\delta_1}^{(2)}\to\C^\times$ by $\chi^{(2,\psi)}_{\delta_1}(X):=\chi^{(1,\psi)}_{\delta_1}(X)\,\psi(\vp g^{-1}h)$. Then:
\begin{enumerate}
\item[(1)] $\chi^{(1,\psi)}_{\delta_1}$ and $\chi^{(2,\psi)}_{\delta_1}$ are well-defined characters of $B_{\delta_1}^{(1)}$ and $B_{\delta_1}^{(2)}$, respectively;
\item[(2)] {\chg we have
\[
\Ind_{B_{\delta_1}}^{B_1}\chi^{\delta_1,\psi}\;\cong\;\Ind_{B_{\delta_1}^{(1)}}^{B_1}\chi^{(1,\psi)}_{\delta_1}\;\oplus\;\Ind_{B_{\delta_1}^{(2)}}^{B_1}\chi^{(2,\psi)}_{\delta_1},
\]
and the two summands are irreducible and non-isomorphic representations of $B_1$.}
\end{enumerate}
\end{proposition}

\begin{proof}
Part (1) follows as in the proofs of {\chg Proposition~\ref{prop:B21} and Lemma~\ref{lem:B2132}}, using $a\equiv e$ on $B_{\delta_1}^{(1)}$ and $e\equiv g$ on $B_{\delta_1}^{(2)}$.

{\chg For (2) we apply Lemma~\ref{lem:split} with $H=B_{\delta_1}$, $H^+=B_{\delta_1}^{(1)}$, $\lambda^+=\chi^{(1,\psi)}_{\delta_1}$, $H'=B_{\delta_1}^{(2)}$ and $\theta(X):=\psi(\vp g^{-1}h)$, so that $\lambda^+|_{H'}\theta=\chi^{(2,\psi)}_{\delta_1}$. Since $h=0$ on $B_{\delta_1}$ and $\ct|_{B_{21}}=\chi^{(21)}$, we have $\lambda^+|_{B_{\delta_1}}=\chi^{\delta_1,\psi}$. Here $[B_{\delta_1}^{(1)}:B_{\delta_1}]=q$, since $h$ only matters modulo $\vp$. The subgroup $B_{\delta_1}^{(2)}$ is the kernel of the surjective homomorphism $B_{\delta_1}^{(1)}\to\Fq^\times$, $X\mapsto eg^{-1}\bmod\vp$, and the elements $\gamma(x):=\diag(1,1,x)\in B_{\delta_1}$, $x\in\cO_2^\times$, represent its cosets; hence $B_{\delta_1}^{(1)}=B_{\delta_1}B_{\delta_1}^{(2)}$. By (1), $\theta$ is a character of $B_{\delta_1}^{(2)}$; it is trivial on $B_{\delta_1}\cap B_{\delta_1}^{(2)}$, where $h=0$, and non-trivial. Finally $\theta^{\gamma(x)}(X)=\psi(\vp g^{-1}hx^{-1})$, so as before $\theta^x\neq\theta$ for every $x\in B_{\delta_1}^{(1)}\smallsetminus B_{\delta_1}^{(2)}$. Lemma~\ref{lem:split} gives
\[
\Ind_{B_{\delta_1}}^{B_{\delta_1}^{(1)}}\chi^{\delta_1,\psi}\cong\chi^{(1,\psi)}_{\delta_1}\oplus\Ind_{B_{\delta_1}^{(2)}}^{B_{\delta_1}^{(1)}}\chi^{(2,\psi)}_{\delta_1},
\]
and inducing to $B_1$ we obtain the isomorphism in (2). By Theorem~\ref{thm:B1}, its left-hand side is the direct sum of two non-isomorphic irreducible representations. Since both summands on the right-hand side are non-zero, they are irreducible and non-isomorphic.}
\end{proof}

Using transitivity of induction, we obtain
\begin{equation}\label{eq:Ichi}
I(\chi)\;\cong\;\Ind_{B_1}^G\ct\;\oplus\;\Ind_{B_{\delta_3}}^G\chi^{\delta_3,\psi}\;\oplus\;\Ind_{B_{\delta_2}}^G\chi^{\delta_2,\psi}\;\oplus\;\Ind_{B_{\delta_1}}^G\chi^{\delta_1,\psi}.
\end{equation}
We now check the irreducibility of the components of $I(\chi)$ in~\eqref{eq:Ichi}.

\begin{lemma}\label{lem:Gchitilde}
Let $\chi_1,\chi_2,\chi_3\in\Oh$ with $\chi_1\sim\chi_2\sim\chi_3$. Then
\[
\dim\End_G\bigl(\Ind_{B_1}^G\ct\bigr)=\begin{cases}
1&\text{if }\chi_1,\chi_2,\chi_3\text{ are pairwise distinct},\\
2&\text{if }\chi_i=\chi_j\neq\chi_k\text{ for }\{i,j,k\}=\{1,2,3\},\\
6&\text{if }\chi_1=\chi_2=\chi_3.
\end{cases}
\]
\end{lemma}

\begin{proof}
Since $B_1$ contains $K$, we have $B_1\backslash G/B_1\cong\Sym_3$. By Mackey theory, $\End_G\bigl(\Ind_{B_1}^G\ct\bigr)\cong\bigoplus_{\delta\in\Sym_3}\Hom_{\delta B_1\delta^{-1}\cap B_1}(\ct^\delta,\ct)$, and it is simple to check that for $\{i,j,k\}=\{1,2,3\}$, we have
\[
\Hom_{\delta B_1\delta^{-1}\cap B_1}\bigl(\ct^\delta,\ct\bigr)\neq0\iff\begin{cases}\chi_i=\chi_j&\text{if }\delta=(ij),\ i\neq j,\\ \chi_1=\chi_2=\chi_3&\text{if }\delta\in\{(123),(132)\}.\end{cases}
\]
Thus, the lemma follows.
\end{proof}

\begin{theorem}\label{thm:Bdelta}
Let $\chi_1,\chi_2,\chi_3\in\Oh$ be pairwise distinct with $\chi_1\sim\chi_2\sim\chi_3$. Then:
\begin{enumerate}
\item[(1)] the representation $\Ind_{B_{\delta_3}}^G\chi^{\delta_3,\psi}$ is irreducible;
\item[(2)] the representation $\Ind_{B_{\delta_2}}^G\chi^{\delta_2,\psi}$ has exactly two irreducible components;
\item[(3)] the representation $\Ind_{B_{\delta_1}}^G\chi^{\delta_1,\psi}$ has exactly three irreducible components.
\end{enumerate}
\end{theorem}

\begin{proof}
We use Mackey theory. For $i\in\{1,2,3\}$ we have
\begin{equation}\label{eq:Wdelta}
\End_G\bigl(\Ind_{B_{\delta_i}}^G\chi^{\delta_i,\psi}\bigr)\cong\Hom_G\Bigl(\Ind_{B_1}^G\Ind_{B_{\delta_i}}^{B_1}\chi^{\delta_i,\psi},\ \Ind_{B_{\delta_i}}^G\chi^{\delta_i,\psi}\Bigr)\cong\bigoplus_{\gamma\in[B_1\backslash G/B_{\delta_i}]}W_{\delta_i,\gamma},
\end{equation}
where
\[
W_{\delta_i,\gamma}:=\Hom_{B_1}\Bigl(\Ind_{B_{\delta_i}}^{B_1}\chi^{\delta_i,\psi},\ \Ind_{\gamma B_{\delta_i}\gamma^{-1}\cap B_1}^{B_1}\bigl(\chi^{\delta_i,\psi}\bigr)^\gamma\Bigr).
\]
Since $K\subseteq B_1$, for each $i$ we have $B_1\backslash G/B_{\delta_i}\cong\overline{B}\backslash\GL_3(\Fq)/\overline{B}_{\delta_i}\cong\Sym_3$, and we take the permutation matrices as representatives. For $\gamma\neq e$, Mackey theory gives
\[
W_{\delta_i,\gamma}\cong\bigoplus_{\beta\in[(\gamma B_{\delta_i}\gamma^{-1}\cap B_1)\backslash B_1/B_{\delta_i}]}\Hom_{\beta B_{\delta_i}\beta^{-1}\cap(\gamma B_{\delta_i}\gamma^{-1}\cap B_1)}\Bigl(\bigl(\chi^{\delta_i,\psi}\bigr)^\beta,\ \bigl(\chi^{\delta_i,\psi}\bigr)^\gamma\Bigr).
\]

(1) Take $i=3$. By Theorem~\ref{thm:B1}, $W_{\delta_3,e}=\End_{B_1}\bigl(\Ind_{B_{\delta_3}}^{B_1}\chi^{\delta_3,\psi}\bigr)\cong\C$. We claim that $W_{\delta_3,\gamma}=0$ for $\gamma\neq e$. If $\gamma\in\{(12),(123),(132)\}$, then $(\gamma B_{\delta_3}\gamma^{-1}\cap B_1)\backslash B_1/B_{\delta_3}$ is a single double coset, and one checks that
\[
W_{\delta_3,\gamma}\cong\Hom_{B_{\delta_3}\cap(\gamma B_{\delta_3}\gamma^{-1}\cap B_1)}\Bigl(\chi^{\delta_3,\psi},\ \bigl(\chi^{\delta_3,\psi}\bigr)^\gamma\Bigr)=0.
\]
Now let $\gamma\in\{(13),(23)\}$. The elements $\beta(x):=\diag(x,1,1)$, $x\in\cO_2^\times$, form an exhaustive set of representatives for $(\gamma B_{\delta_3}\gamma^{-1}\cap B_1)\backslash B_1/B_{\delta_3}$. For $\gamma=(13)$ and $\beta=\beta(x)$, the group $\beta B_{\delta_3}\beta^{-1}\cap(\gamma B_{\delta_3}\gamma^{-1}\cap B_1)$ contains $\{I_3+\vp cE_{13}:c\in\cO_2\}$. The character $(\chi^{\delta_3,\psi})^\beta$ is trivial on this subgroup, while $(\chi^{\delta_3,\psi})^\gamma(I_3+\vp cE_{13})=\chi^{\delta_3,\psi}(I_3+\vp cE_{31})=\psi(\vp c)$ is not; thus $W_{\delta_3,(13)}=0$. For $\gamma=(23)$ and $\beta=\beta(x)$, the group $\beta B_{\delta_3}\beta^{-1}\cap(\gamma B_{\delta_3}\gamma^{-1}\cap B_1)$ contains $\{I_3+\vp cE_{21}:c\in\cO_2\}$, on which $(\chi^{\delta_3,\psi})^\beta$ is trivial while $(\chi^{\delta_3,\psi})^\gamma$ is not. Thus $W_{\delta_3,(23)}=0$, and~\eqref{eq:Wdelta} gives $\End_G\bigl(\Ind_{B_{\delta_3}}^G\chi^{\delta_3,\psi}\bigr)\cong\C$.

(2) Take $i=2$. 
By Theorem~\ref{thm:B1}, $W_{\delta_2,e}=\End_{B_1}\bigl(\Ind_{B_{\delta_2}}^{B_1}\chi^{\delta_2,\psi}\bigr)\cong\C$. 
We claim that $W_{\delta_2,\gamma}=0$ for $\gamma\notin\{e,(12)\}$ and $W_{\delta_2,(12)}\cong\C$. For $\gamma\in\{(13),(123),(132)\}$ one checks, as in (1), that $W_{\delta_2,\gamma}\cong\Hom_{B_{\delta_2}\cap(\gamma B_{\delta_2}\gamma^{-1}\cap B_1)}\bigl(\chi^{\delta_2,\psi},(\chi^{\delta_2,\psi})^\gamma\bigr)=0$.

 For $\gamma=(23)$, the elements $\beta(x):=\diag(1,x,1)$, $x\in\cO_2^\times$, form an exhaustive set of representatives for $(\gamma B_{\delta_2}\gamma^{-1}\cap B_1)\backslash B_1/B_{\delta_2}$.
For $\beta=\beta(x)$, the group $\beta B_{\delta_2}\beta^{-1}\cap(\gamma B_{\delta_2}\gamma^{-1}\cap B_1)$ contains $\{I_3+\vp cE_{23}:c\in\cO_2\}$. The character $(\chi^{\delta_2,\psi})^\beta$ is trivial on it, while $(\chi^{\delta_2,\psi})^\gamma$ is not, so $W_{\delta_2,(23)}=0$. Finally, for $\gamma=(12)$ the double coset space $(\gamma B_{\delta_2}\gamma^{-1}\cap B_1)\backslash B_1/B_{\delta_2}$ is a single point and
\[
W_{\delta_2,(12)}\cong\Hom_{B_{\delta_2}\cap(\gamma B_{\delta_2}\gamma^{-1}\cap B_1)}\Bigl(\chi^{\delta_2,\psi},\ \bigl(\chi^{\delta_2,\psi}\bigr)^\gamma\Bigr),
\]
where
\[
B_{\delta_2}\cap\bigl(\gamma B_{\delta_2}\gamma^{-1}\cap B_1\bigr)=\left\{A:=\begin{pmatrix}a&\vp b&c\\ \vp d&e&f\\ 0&0&g\end{pmatrix}:\begin{array}{l}a,e,g\in\cO_2^\times,\ a\equiv e\equiv g\bmod\vp,\\ b,c,d,f\in\cO_2\end{array}\right\}.
\]
On this group $\chi^{\delta_2,\psi}(A)=\chi_1(a)\chi_2(e)\chi_3(g)$ and $(\chi^{\delta_2,\psi})^\gamma(A)=\chi_1(e)\chi_2(a)\chi_3(g)$; these agree because $\chi_1\sim\chi_2$ and $ae^{-1}\in1+\vp\cO_2$. Thus $W_{\delta_2,(12)}\cong\C$, and~\eqref{eq:Wdelta} gives $\dim\End_G\bigl(\Ind_{B_{\delta_2}}^G\chi^{\delta_2,\psi}\bigr)=2$.

(3) Take $i=1$.
By Theorem~\ref{thm:B1}, $\dim W_{\delta_1,e}=\dim\End_{B_1}\bigl(\Ind_{B_{\delta_1}}^{B_1}\chi^{\delta_1,\psi}\bigr)=2$. 
We claim that $W_{\delta_1,\gamma}=0$ for $\gamma\notin\{e,(23)\}$ and $W_{\delta_1,(23)}\cong\C$. 
For $\gamma\in\{(123),(132)\}$ one checks that $W_{\delta_1,\gamma}\cong\Hom_{B_{\delta_1}\cap(\gamma B_{\delta_1}\gamma^{-1}\cap B_1)}\bigl(\chi^{\delta_1,\psi},(\chi^{\delta_1,\psi})^\gamma\bigr)=0$.
Now let $\gamma=(12)$. 
An exhaustive set of representatives for $(\gamma B_{\delta_1}\gamma^{-1}\cap B_1)\backslash B_1/B_{\delta_1}$ is given by $\beta(x):=\diag(x,1,1)$, $x\in\cO_2^\times$, together with $I_3+\vp E_{31}$ and $I_3+\vp E_{32}$. 
For $\beta=\beta(x)$, the group $\beta B_{\delta_1}\beta^{-1}\cap(\gamma B_{\delta_1}\gamma^{-1}\cap B_1)$ contains $\{I_3+\vp cE_{21}:c\in\cO_2\}$; the character $(\chi^{\delta_1,\psi})^\gamma$ is trivial on it, while $(\chi^{\delta_1,\psi})^\beta(I_3+\vp cE_{21})=\psi(\vp cx)$ is not.
Similarly, for $\beta\in\{I_3+\vp E_{31},I_3+\vp E_{32}\}$ the group $\beta B_{\delta_1}\beta^{-1}\cap(\gamma B_{\delta_1}\gamma^{-1}\cap B_1)$ contains $\{I_3+\vp cE_{21}:c\in\cO_2\}$, on which $(\chi^{\delta_1,\psi})^\gamma$ is trivial while $(\chi^{\delta_1,\psi})^\beta$ is not.
Thus $W_{\delta_1,(12)}=0$, and a similar argument shows that $W_{\delta_1,(13)}=0$. Finally, for $\gamma=(23)$ the double coset space $(\gamma B_{\delta_1}\gamma^{-1}\cap B_1)\backslash B_1/B_{\delta_1}$ is a single point, so
\[
W_{\delta_1,(23)}\cong\Hom_{B_{\delta_1}\cap(\gamma B_{\delta_1}\gamma^{-1}\cap B_1)}\Bigl(\chi^{\delta_1,\psi},\ \bigl(\chi^{\delta_1,\psi}\bigr)^\gamma\Bigr),
\]
where
\[
B_{\delta_1}\cap\bigl(\gamma B_{\delta_1}\gamma^{-1}\cap B_1\bigr)=\left\{\begin{pmatrix}a&b&c\\ 0&e&0\\ 0&0&g\end{pmatrix}:a,e,g\in\cO_2^\times,\ a\equiv e\equiv g\bmod\vp,\ b,c\in\cO_2\right\}.
\]
On this group $\chi^{\delta_1,\psi}$ and $(\chi^{\delta_1,\psi})^\gamma$ are $\chi_1(a)\chi_2(e)\chi_3(g)$ and $\chi_1(a)\chi_2(g)\chi_3(e)$, which agree since $\chi_2\sim\chi_3$.
Hence $W_{\delta_1,(23)}\cong\C$, and~\eqref{eq:Wdelta} gives $\dim\End_G\bigl(\Ind_{B_{\delta_1}}^G\chi^{\delta_1,\psi}\bigr)=3$.

Thus, the theorem follows.
\qedhere
\end{proof}

\subsection{A description of \texorpdfstring{$\Ind_{B_{\delta_2}}^G\chi^{\delta_2,\psi}$}{Ind(B delta2,G)}}
In this subsection we describe the two irreducible components of $\Ind_{B_{\delta_2}}^G\chi^{\delta_2,\psi}$. We introduce the following two subgroups of $B_1$:
\[
B':=\left\{X=\begin{pmatrix}a&b&c\\ \vp d&e&f\\ \vp k&\vp h&g\end{pmatrix}\in B_1:a-g\in\vp\cO_2\right\},\qquad B'':=\bigl\{X\in B':e-g\in\vp\cO_2\bigr\}.
\]
Note that $B'$ coincides with $B_{\delta_3}$ as a group; we keep the separate notation because it carries a different character in the Lemma below. Note also that $B'$ is the inverse image in $G$ of the centralizer of $E_{13}$ in $\GL_3(\Fq)$.

The following lemma is elementary.

\begin{lemma}\label{lem:chiprime}
Let $\chi_1,\chi_2,\chi_3\in\Oh$ with $\chi_1\sim\chi_2\sim\chi_3$, and let $\psi$ be as above. Define $\chi^{\prime,\psi}\colon B'\to\C^\times$ and $\chi^{\prime \prime, \psi}\colon B''\to\C^\times$ by
\begin{align*}
\chi^{\prime, \psi}(X)
&:=\chi_2(a)\chi_1(e)\chi_3(g)\cdot\chi_1\bigl(1-\vp a^{-1}e^{-1}b(d-g^{-1}fk)-\vp e^{-1}g^{-1}fh\bigr)\\
&\qquad\cdot\chi_1\bigl(1-\vp a^{-1}g^{-1}ck\bigr)\cdot\psi(\vp g^{-1}k),\\
\chi^{\prime \prime, \psi}(X)&:=\chi^{\prime, \psi}(X)\,\psi(\vp e^{-1}b).
\end{align*}
Equivalently, $\chi^{\prime, \psi}(X)=\chi_1(\det X)\,(\chi_1^{-1}\chi_2)(a)\,(\chi_1^{-1}\chi_3)(g)\,\psi(\vp g^{-1}k)$. 
Then $\chi^{\prime, \psi}$ and $\chi^{\prime \prime, \psi}$ are well-defined characters of $B'$ and $B''$, respectively.
\end{lemma}

\begin{proof}
The equivalent form follows from
\[
\det X=aeg\bigl(1-\vp a^{-1}e^{-1}b(d-g^{-1}fk)-\vp e^{-1}g^{-1}fh-\vp a^{-1}g^{-1}ck\bigr).
\]
The map $X\mapsto\chi_1(\det X)(\chi_1^{-1}\chi_2)(a)(\chi_1^{-1}\chi_3)(g)$ is a character of $B_1$, as in Lemma~\ref{lem:chitilde}. 
The map $X\mapsto\psi(\vp g^{-1}k)$ is a character of $B'$ because $a\equiv g\bmod\vp$ there {\chg (see the proof of Lemma~\ref{lem:tildetau})}. Finally, for $X,X'\in B''$ we have $\vp(XX')_{12}=\vp(ab'+be')$ and $(XX')_{22}\equiv ee'\bmod\vp$; since $a\equiv e\bmod\vp$ on $B''$, the map $X\mapsto\psi(\vp e^{-1}b)$ is a character of $B''$.
\end{proof}

\begin{proposition}\label{prop:BprimeBsecond}
Let $\chi_1,\chi_2,\chi_3\in\Oh$ with $\chi_1\sim\chi_2\sim\chi_3$. Then
\[
\Ind_{B_{\delta_2}}^G\chi^{\delta_2,\psi}\;\cong\;\Ind_{B'}^G\chi^{\prime, \psi}\;\oplus\;\Ind_{B''}^G\chi^{\prime \prime, \psi}.
\]
The two summands are irreducible, of dimensions $(q+1)(q^3-1)$ and $(q^2-1)(q^3-1)$.
\end{proposition}

\begin{proof}
Write $\tau_1:=\Ind_{B'}^G\chi^{\prime, \psi}$ and $\tau_2:=\Ind_{B''}^G\chi^{\prime\prime,\psi}$. Since $\dim\Ind_{B_{\delta_2}}^G\chi^{\delta_2,\psi}=q(q+1)(q^3-1)=\dim\tau_1+\dim\tau_2$ and $\tau_1\not\cong\tau_2$, the proposition follows from the following two statements:
\begin{enumerate}
\item[(1)] $\End_G(\tau_i)\cong\C$ for $i=1,2$;
\item[(2)] $\Hom_G\bigl(\Ind_{B_{\delta_2}}^G\chi^{\delta_2,\psi},\tau_i\bigr)\neq0$ for $i=1,2$.
\end{enumerate}
Since the proofs for $i=1$ and $i=2$ are analogous, we give the proof for $i=1$; for $i=2$ one uses the representatives $\{\diag(1,x,y):x,y\in\cO_2^\times\}$ of $B''\backslash B_1/B''$.

\emph{Proof of (1).} By Mackey theory,
\[
\End_G(\tau_1)\cong\bigoplus_{\delta\in[B_1\backslash G/B']}W_\delta,\qquad W_\delta:=\Hom_{B_1}\Bigl(\Ind_{B'}^{B_1}\chi^{\prime, \psi},\ \Ind_{\delta B'\delta^{-1}\cap B_1}^{B_1}(\chi^{\prime, \psi})^\delta\Bigr),
\]
and $B_1\backslash G/B'\cong\Sym_3$.
We claim that $W_e\cong\C$ and $W_\delta=0$ for $\delta\neq e$. 
Suppose first that $\delta=e$. 
Then
\[
W_e\cong\Hom_{B'}\Bigl(\Res_{B'}\Ind_{B'}^{B_1}\chi^{\prime, \psi},\chi^{\prime, \psi}\Bigr)\cong\bigoplus_{\gamma\in[B'\backslash B_1/B']}\Hom_{\gamma B'\gamma^{-1}\cap B'}\bigl((\chi^{\prime, \psi})^\gamma,\chi^{\prime, \psi}\bigr).
\]
The subgroup $B'$ is normal in $B_1$, a set of representatives for $B'\backslash B_1/B'$ is $\{\gamma(x):=\diag(1,1,x):x\in\cO_2^\times\}$, and $B'\gamma(x)B'=B'\gamma(y)B'$ if and only if $x\equiv y\bmod\vp$. 
For $\gamma=\gamma(x)$ and $X\in B'$,
\[
\chi^{\prime, \psi}(A)\,(\chi^{\prime, \psi})^\gamma(A)^{-1}=\psi\bigl(\vp(1-x^{-1})g^{-1}k\bigr)\implies(\chi^{\prime, \psi})^\gamma=\chi^{\prime,\psi}\iff x\in1+\vp\cO_2.
\]
This proves $W_e\cong\C$.

Now let $\delta\neq e$. By Mackey theory,
\[
W_\delta\cong\bigoplus_{\beta\in[B'\backslash B_1/(\delta B'\delta^{-1}\cap B_1)]}\Hom_{\beta\delta B'(\beta\delta)^{-1}\cap B'}\bigl(\chi^{\prime, \psi},(\chi^{\prime, \psi})^{\beta\delta}\bigr).
\]
The group $\beta\delta B'(\beta\delta)^{-1}\cap B'$ contains $K$.
Moreover for $Y\in M_3(\cO_2)$, $\chi^{\prime, \psi}(I_3+\vp Y)=\chi_1(1+\vp\tr Y)\,\psi(\vp Y_{31})$, while $(\chi^{\prime, \psi})^{\beta\delta}(I_3+\vp Y)=\chi_1(1+\vp\tr Y)\,\psi\bigl(\vp((\beta\delta)^{-1}Y\beta\delta)_{31}\bigr)$. 
Since $Y_{31}=\tr(E_{13}Y)$, these two characters of $K$ agree if and only if $\beta\delta$ modulo $\vp$ commutes with $E_{13}$, that is, if and only if $\beta\delta\in B'$. This is because $B'$ is the inverse image in $G$ of the centralizer of $E_{13}$ in $\GL_3(\Fq)$. This is not the case, since $\beta\in B_1$ and $\delta\notin B_1$. Hence $W_\delta=0$, which proves (1).

\emph{Proof of (2).} By Mackey theory,
\[
\Hom_G\bigl(\Ind_{B_{\delta_2}}^G\chi^{\delta_2,\psi},\tau_1\bigr)\cong\Hom_{B_1}\Bigl(\Ind_{B_{\delta_2}}^{B_1}\chi^{\delta_2,\psi},\ \Res_{B_1}\Ind_{B'}^G\chi^{\prime, \psi}\Bigr)\cong\bigoplus_{\gamma\in[B_1\backslash G/B']}W_{\delta_2,\gamma},
\]
where $W_{\delta_2,\gamma}:=\Hom_{B_1}\bigl(\Ind_{B_{\delta_2}}^{B_1}\chi^{\delta_2,\psi},\ \Ind_{\gamma B'\gamma^{-1}\cap B_1}^{B_1}(\chi^{\prime, \psi})^\gamma\bigr)$. 
Let $\gamma=(12)$, regarded as the permutation matrix $P_{(12)}$.
By Mackey theory,
\[
W_{\delta_2,\gamma}\cong\bigoplus_{\beta\in[B_{\delta_2}\backslash B_1/(\gamma B'\gamma^{-1}\cap B_1)]}\Hom_{\beta\gamma B'(\beta\gamma)^{-1}\cap B_{\delta_2}}\bigl(\chi^{\delta_2,\psi},(\chi^{\prime, \psi})^{\beta\gamma}\bigr).
\]
For $\beta=I_3$,
\[
\gamma B'\gamma^{-1}\cap B_{\delta_2}=\left\{A=\begin{pmatrix}a&\vp b&c\\ \vp d&e&f\\ 0&\vp h&g\end{pmatrix}:\begin{array}{l}a,e,g\in\cO_2^\times,\ e\equiv g\bmod\vp,\\ b,c,d,f,h\in\cO_2\end{array}\right\},
\]
and
\[
(\chi^{\prime, \psi})^\gamma(A)=\chi^{\prime, \psi}(\gamma^{-1}A\gamma)=\chi^{\prime, \psi}\begin{pmatrix}
e&\vp d&f\\
\vp b&a&c\\
\vp h&0&g\end{pmatrix}.
\]
Using the second description of $\chi'$ in Lemma~\ref{lem:chiprime}, one checks that $\chi^{\delta_2,\psi}(A)=(\chi^{\prime, \psi})^\gamma(A)$ for all such $A$; in fact both are equal to
\[
\chi_1(\det A)\,(\chi_2\chi_3\chi_1^{-2})(e)\,\psi(\vp g^{-1}h).
\]
Hence $W_{\delta_2,(12)}\neq0$, which proves (2).
\end{proof}

\begin{lemma}\label{lem:new-chars}
Let $\chi_1,\chi_2,\chi_3\in\Oh$ with $\chi_1\sim\chi_2\sim\chi_3$, and write $X\in B'$ as
\[
X=\begin{pmatrix} a&b&c\\ \varpi d&e&f\\ \varpi k&\varpi h&g\end{pmatrix},\qquad a\equiv g \pmod{\varpi}.
\]
Define
\begin{align*}
\chi^{*}_{\psi}(X)&:=\chi_1(\det X)\,(\chi_1^{-1}\chi_2)(a)\,(\chi_1^{-1}\chi_3)(e)\,\psi(\varpi g^{-1}k) && (X\in B'),\\
\chi^{**}_{\psi}(X)&:=\chi^{\prime,\psi}(X)\,\psi(\varpi g^{-1}f) && (X\in B'').
\end{align*}
Then $\chi^{*}_{\psi}$ is a character of $B'$ and $\chi^{**}_{\psi}$ is a character of $B''$.
\end{lemma}

\begin{proof}
As in Lemma~5.1, $X\mapsto\chi_1(\det X)(\chi_1^{-1}\chi_2)(a)(\chi_1^{-1}\chi_3)(e)$ is a character of $B_1$, and $X\mapsto\psi(\varpi g^{-1}k)$ is a character of $B'$ by the proof of Lemma~5.12.
For $X,X'\in B''$ we have $(XX')_{23}=\varpi dc'+ef'+fg'\equiv ef'+fg'$ and $(XX')_{33}\equiv gg' \pmod{\varpi}$. Since $e\equiv g\pmod\varpi$ on $B''$,
\[
\frac{ef'+fg'}{gg'}\equiv\frac{f'}{g'}+\frac{f}{g}\pmod{\varpi},
\]
so $X\mapsto\psi(\varpi g^{-1}f)$ is a character of $B''$, and hence so is $\chi^{**}_{\psi}$.
\end{proof}

Note that $\chi^{*}_{\psi}$ is obtained from $\chi^{\prime,\psi}$ by replacing $(\chi_1^{-1}\chi_3)(g)$ with $(\chi_1^{-1}\chi_3)(e)$, and $\chi^{**}_{\psi}$ is obtained from $\chi^{\prime\prime,\psi}$ by replacing $\psi(\varpi e^{-1}b)$ with $\psi(\varpi g^{-1}f)$.

\begin{proposition}\label{prop:delta1-split}
Let $\chi_1,\chi_2,\chi_3\in\Oh$ with $\chi_1\sim\chi_2\sim\chi_3$. 
Then
\[
\Ind_{B_{\delta_1}^{(1)}}^G\chi^{(1,\psi)}_{\delta_1}\cong\Ind_{B'}^G\chi^{*}_{\psi}\oplus\Ind_{B''}^G\chi^{**}_{\psi},
\]
where the two summands are irreducible and of dimensions $(q+1)(q^3-1)$ and $(q^2-1)(q^3-1)$.
\end{proposition}

\begin{proof}
The irreducibility of $\Ind_{B'}^{G}\chi^*_{\psi}$ and $\Ind_{B'}^{G}\chi^{**}_{\psi}$ is same as irreducibility of the representations \(\Ind_{B'}^G\chi^{\prime, \psi}\) and \(Ind_{B''}^G\chi^{\prime \prime, \psi}\) in Proposition \ref{prop:BprimeBsecond}, so we skip the proof.

Now, we prove that both the mentioned representations occur in \(\Ind_{B_{\delta_1}^{(1)}}^G\chi^{(1,\psi)}_{\delta_1}\).
Recall that \[
B_{\delta_1}^{(1)} :=\left\{X=\begin{pmatrix}a&b&c\\ \vp d&e&f\\ 0&\vp h&g\end{pmatrix}\in B_{21,32}:a-e\in\vp\cO_2\right\}.\]
Let $s=P_{(23)}$.
Then
\(
B^{(1)}_{\delta_1}\cap sB's^{-1}=\{X\in B^{(1)}_{\delta_1}: f\in\varpi\mathcal O_2\},
\)
and for such $X$,
\[
Z:=s^{-1}Xs=\begin{pmatrix} a&c&b\\ 0&g&\varpi h \\ \varpi d& f&e\end{pmatrix}\in B'.
\]
Hence
\[
(\chi^{*}_{\psi})^{s}(X)=\chi^{*}_{\psi}(Z)=\chi_1(\det X)(\chi_1^{-1}\chi_2)(a)(\chi_1^{-1}\chi_3)(g)\,\psi(\varpi e^{-1}d)=\chi^{(1,\psi)}_{\delta_1}(X),
\]
since $a\equiv e\pmod\varpi$ on $B^{(1)}_{\delta_1}$. 
So the term of $s$ in Mackey's formula for
\(
\Hom_G\big(\Ind_{B^{(1)}_{\delta_1}}^G\chi^{(1,\psi)}_{\delta_1},\Ind_{B'}^G\chi^{*}_{\psi}\big)\neq 0.
\)

Similarly, $B^{(1)}_{\delta_1}\cap sB''s^{-1}=\{X\in B^{(1)}_{\delta_1}: f\in\varpi\mathcal O_2,\ a\equiv g\}$.
For such $X$, the $(2,3)$-entry of $Z$ is $f\in\varpi\mathcal O_2$ and
\[
(\chi^{**}_{\psi})^{s}(X)=\chi^{\prime,\psi}(Z)=\chi_1(\det X)(\chi_1^{-1}\chi_2)(a)(\chi_1^{-1}\chi_3)(e)\,\psi(\varpi e^{-1}d)=\chi^{(1,\psi)}_{\delta_1}(X),
\]
since $a\equiv e\equiv g\pmod\varpi$.
Hence using Mackey theory, we get
\(
\Hom_G\big(\Ind_{B^{(1)}_{\delta_1}}^G\chi^{(1,\psi)}_{\delta_1},\Ind_{B'}^G\chi^{**}_{\psi}\big)\neq 0.
\)
Since
\[
(q+1)(q^3-1)+(q^2-1)(q^3-1)=q(q+1)(q^3-1)=\dim\Ind_{B^{(1)}_{\delta_1}}^G\chi^{(1,\psi)}_{\delta_1},
\]
the proposition follows.
\end{proof}

\begin{theorem}\label{thm:D7}
Let $\chi_1,\chi_2,\chi_3\in\Oh$ be pairwise distinct with $\chi_1\sim\chi_2\sim\chi_3$. Then
\begin{align*}
I(\chi)& \cong\Ind_{B_1}^G\ct\oplus\Ind_{B_{\delta_3}}^G\chi^{\delta_3,\psi}\oplus\Ind_{B'}^G\chi^{\prime,\psi}\oplus\Ind_{B''}^G\chi^{\prime\prime,\psi}
\oplus\Ind_{B'}^G\chi^{*}_{\psi}\oplus\Ind_{B''}^G\chi^{**}_{\psi}\oplus\Ind_{B^{(2)}_{\delta_1}}^G\chi^{(2,\psi)}_{\delta_1},
\end{align*}
and the 7 summands are irreducible and pairwise non-isomorphic.
\end{theorem}

\begin{proof}
By Proposition~5.13, Proposition~5.9 and Proposition~\ref{prop:delta1-split}, $I(\chi)$ is the direct sum of the seven non-zero representations shown. Since $D(\chi)=7$ by Proposition~3.4, they are irreducible and pairwise non-isomorphic.
\end{proof}

\begin{table}[ht]
\caption{The seven irreducible constituents of $I(\chi)$ when $D(\chi)=7$.}\label{tab:D7}
\renewcommand{\arraystretch}{1.3}
\begin{tabular}{lll}
\hline
Irreducible representation & Dimension\\
\hline
$\Ind_{B_1}^G\ct$ & $(q+1)(q^2+q+1)$\\
$\Ind_{B_{\delta_3}}^G\chi^{\delta_3,\psi}$ & $(q+1)(q^3-1)$\\
$\Ind_{B'}^G\chi^{\prime,\psi}$ & $(q+1)(q^3-1)$\\
$\Ind_{B'}^G\chi^{*}_{\psi}$ & $(q+1)(q^3-1)$\\
$\Ind_{B''}^G\chi^{\prime\prime,\psi}$ & $(q^2-1)(q^3-1)$\\
$\Ind_{B''}^G\chi^{**}_{\psi}$ & $(q^2-1)(q^3-1)$\\
$\Ind_{B^{(2)}_{\delta_1}}^G\chi^{(2,\psi)}_{\delta_1}$ & $q(q^2-1)(q^3-1)$\\
\end{tabular}
\end{table}

\section{The cases \texorpdfstring{$D(\chi)\in\{10,18\}$}{D(chi) in \{10,18\}}}\label{sec:D1018}

We keep the assumption $\chi_1\sim\chi_2\sim\chi_3$. By Theorem~\ref{thm:D7}, the representations occurring in $I(\chi)$ when $D(\chi)=7$, together with their dimensions, are listed in Table~\ref{tab:D7}.

\begin{proposition}\label{prop:homs}
Let $\chi_1,\chi_2,\chi_3\in\Oh$ with $\chi_1\sim\chi_2\sim\chi_3$. Then:
\begin{enumerate}
\item[(1)] $\Hom_G\bigl(\Ind_{B_{\delta_3}}^G\chi^{\delta_3,\psi},\ \Ind_{B_{\delta_2}}^G\chi^{\delta_2,\psi}\bigr)\cong\C$ if $\chi_1=\chi_2$, and is $0$ otherwise;
\item[(2)] $\Hom_G\bigl(\Ind_{B_{\delta_3}}^G\chi^{\delta_3,\psi},\ \Ind_{B_{\delta_1}}^G\chi^{\delta_1,\psi}\bigr)\cong\C$ if $\chi_2=\chi_3$, and is $0$ otherwise;
\item[(3)] $\Hom_G\bigl(\Ind_{B_{\delta_2}}^G\chi^{\delta_2,\psi},\ \Ind_{B_{\delta_1}}^G\chi^{\delta_1,\psi}\bigr)\cong\C$ if $\chi_1=\chi_3$, and is $0$ otherwise.
\end{enumerate}
\end{proposition}

\begin{proof}
By Frobenius reciprocity and Mackey theory,
\begin{align*}
\Hom_G\Bigl(\Ind_{B_{\delta_3}}^G\chi^{\delta_3,\psi},\ \Ind_{B_{\delta_2}}^G\chi^{\delta_2,\psi}\Bigr)&\cong\Hom_{B_1}\Bigl(\Ind_{B_{\delta_3}}^{B_1}\chi^{\delta_3,\psi},\ \Res_{B_1}\Ind_{B_{\delta_2}}^G\chi^{\delta_2,\psi}\Bigr)\\
&\cong\bigoplus_{\gamma\in[B_1\backslash G/B_{\delta_2}]}\Hom_{B_1}\Bigl(\Ind_{B_{\delta_3}}^{B_1}\chi^{\delta_3,\psi},\ \Ind_{\gamma B_{\delta_2}\gamma^{-1}\cap B_1}^{B_1}\bigl(\chi^{\delta_2,\psi}\bigr)^\gamma\Bigr).
\end{align*}
Since $K\subseteq B_1$, we have $B_1\backslash G/B_{\delta_2}\cong\overline{B}\backslash\GL_3(\Fq)/\overline{B}_{\delta_2}\cong\Sym_3=\{e,(12),(13),(23),(123),(132)\}$. Moreover, it can be checked that
\[
\dim\Hom_{B_1}\Bigl(\Ind_{B_{\delta_3}}^{B_1}\chi^{\delta_3,\psi},\ \Ind_{\gamma B_{\delta_2}\gamma^{-1}\cap B_1}^{B_1}\bigl(\chi^{\delta_2,\psi}\bigr)^\gamma\Bigr)=\begin{cases}1&\text{if }\gamma=e\text{ and }\chi_1=\chi_2,\\ 0&\text{otherwise.}\end{cases}
\]
This proves (1). Similarly, it can be checked that
\[
\dim\Hom_{B_1}\Bigl(\Ind_{B_{\delta_3}}^{B_1}\chi^{\delta_3,\psi},\ \Ind_{\gamma B_{\delta_1}\gamma^{-1}\cap B_1}^{B_1}\bigl(\chi^{\delta_1,\psi}\bigr)^\gamma\Bigr)=\begin{cases}1&\text{if }\gamma=(23)\text{ and }\chi_2=\chi_3,\\ 0&\text{otherwise,}\end{cases}
\]
which proves (2). Finally, it can be checked that
\[
\dim\Hom_{B_1}\Bigl(\Ind_{B_{\delta_2}}^{B_1}\chi^{\delta_2,\psi},\ \Ind_{\gamma B_{\delta_1}\gamma^{-1}\cap B_1}^{B_1}\bigl(\chi^{\delta_1,\psi}\bigr)^\gamma\Bigr)=\begin{cases}1&\text{if }\gamma=(123)\text{ and }\chi_1=\chi_3,\\ 0&\text{otherwise,}\end{cases}
\]
which proves (3).
\end{proof}

Let $P_{2,1}^1:=P_{2,1}\cdot K$ be the inverse image in $G$ of the parabolic subgroup $\overline{P}_{2,1}$ of $\GL_3(\Fq)$; thus $B_1\subseteq P_{2,1}^1$ and $[P_{2,1}^1:B_1]=q+1$. As before, $\St_2$ denotes the Steinberg representation of $\GL_2(\Fq)$, regarded as a representation of $P_{2,1}^1$ via $P_{2,1}^1\twoheadrightarrow\overline{P}_{2,1}\twoheadrightarrow\GL_2(\Fq)$.

\begin{lemma}\label{lem:P21}
Let $\chi_1,\chi_2,\chi_3\in\Oh$ satisfy $\chi_1\sim\chi_2\sim\chi_3$ and $\chi_1=\chi_2\neq\chi_3$. Then $\ct$ extends to the character of $P_{2,1}^1$, again denoted $\ct$, given by $\ct(x)=\chi_1(\det x)\,(\chi_1^{-1}\chi_3)(x_{33})$, and
\[
\Ind_{B_1}^G\ct\;\cong\;\Ind_{P_{2,1}^1}^G\ct\;\oplus\;\Ind_{P_{2,1}^1}^G\bigl(\ct\otimes\St_2\bigr).
\]
The two summands are irreducible and of dimensions $q^2+q+1$ and $q(q^2+q+1)$.
\end{lemma}

\begin{proof}
The map $x\mapsto x_{33}\bmod\vp$ is a homomorphism $P_{2,1}^1\to\Fq^\times$, and $\chi_1^{-1}\chi_3$ is trivial on $1+\vp\cO_2$; hence the formula defines a character of $P_{2,1}^1$, which restricts to $\ct$ on $B_1$ by Lemma~\ref{lem:chitilde} (since $\chi_1=\chi_2$). Therefore $\Ind_{B_1}^{P_{2,1}^1}\ct\cong\ct\otimes\Ind_{B_1}^{P_{2,1}^1}\one$. Since $K\subseteq B_1$, the representation $\Ind_{B_1}^{P_{2,1}^1}\one$ is inflated from $\Ind_{\overline{B}}^{\overline{P}_{2,1}}\one$. Writing $\overline{P}_{2,1}=(\GL_2(\Fq)\times\Fq^\times)\ltimes N$ and $\overline{B}=(\overline{B}_2\times\Fq^\times)\ltimes N$, we obtain $\Ind_{\overline{B}}^{\overline{P}_{2,1}}\one\cong\bigl(\Ind_{\overline{B}_2}^{\GL_2(\Fq)}\one\bigr)\boxtimes\one\cong\one\oplus(\St_2\boxtimes\one)$. Hence $\Ind_{B_1}^{P_{2,1}^1}\ct\cong\ct\oplus(\ct\otimes\St_2)$, and the decomposition follows by transitivity of induction. By Lemma~\ref{lem:Gchitilde}, $\dim\End_G\bigl(\Ind_{B_1}^G\ct\bigr)=2$, so the two summands are irreducible and non-isomorphic.
\end{proof}

When $D(\chi)\in\{10,18\}$, the only new irreducible representations in $I(\chi)$ appears frome the summand $\Ind_{B_1}^G\ct$, (Lemma~5.10), and some of the three characters $\chi^{\delta_3,\psi}$, $\chi^{\prime,\psi}$, $\chi^{*}_{\psi}$ of $B'=B_{\delta_3}$ coincide.
The following lemma explain this phenomenon.

\begin{lemma}\label{lem:compare}
Let $\chi_1,\chi_2,\chi_3\in\Oh$ with $\chi_1\sim\chi_2\sim\chi_3$.
On $B'=B_{\delta_3}$ we have $\chi^{\delta_3,\psi}=\chi^{\prime,\psi}$ if and only if $\chi_1=\chi_2$; $\chi^{\delta_3,\psi}=\chi^{*}_{\psi}$ if and only if $\chi_2=\chi_3$; and $\chi^{\prime,\psi}=\chi^{*}_{\psi}$ if and only if $\chi_1=\chi_3$.
\end{lemma}

\begin{proof}
For $X\in B'$ written as in Lemma~5.14, the definitions before Lemma~5.5, Lemma~5.12 and Lemma~5.14 give
\begin{align*}
\chi^{\delta_3,\psi}(X)&=\chi_1(\det X)\,(\chi_1^{-1}\chi_2)(e)\,(\chi_1^{-1}\chi_3)(g)\,\psi(\varpi g^{-1}k),\\
\chi^{\prime,\psi}(X)&=\chi_1(\det X)\,(\chi_1^{-1}\chi_2)(a)\,(\chi_1^{-1}\chi_3)(g)\,\psi(\varpi g^{-1}k),\\
\chi^{*}_{\psi}(X)&=\chi_1(\det X)\,(\chi_1^{-1}\chi_2)(a)\,(\chi_1^{-1}\chi_3)(e)\,\psi(\varpi g^{-1}k).
\end{align*}
Since $a\equiv g\pmod\varpi$ on $B'$, we have
\begin{gather*}
\chi^{\delta_3,\psi}(\chi^{\prime,\psi})^{-1}=(\chi_1^{-1}\chi_2)(ea^{-1}),\qquad
\chi^{\delta_3,\psi}(\chi^{*}_{\psi})^{-1}=(\chi_3^{-1}\chi_2)(ea^{-1}),\\
\chi^{\prime,\psi}(\chi^{*}_{\psi})^{-1}=(\chi_1^{-1}\chi_3)(ae^{-1}).
\end{gather*} 
Hence the lemma follows.
\end{proof}

\begin{theorem}\label{thm:D10}
Let $\chi_1,\chi_2,\chi_3\in\Oh$ with $\chi_1=\chi_2\neq\chi_3$ and $\chi_1\sim\chi_3$. Then $\chi^{\delta_3,\psi}=\chi^{\prime,\psi}\neq\chi^{*}_{\psi}$, and
\begin{align*}
I(\chi)& \cong\Ind_{P^1_{2,1}}^G\ct\oplus\Ind_{P^1_{2,1}}^G(\ct\otimes\mathrm{St}_2)\oplus2\cdot\Ind_{B_{\delta_3}}^G\chi^{\delta_3,\psi}\oplus\Ind_{B'}^G\chi^{*}_{\psi}
\oplus\Ind_{B''}^G\chi^{\prime\prime,\psi}\oplus\Ind_{B''}^G\chi^{**}_{\psi}
\\
& \oplus\Ind_{B^{(2)}_{\delta_1}}^G\chi^{(2,\psi)}_{\delta_1}.
\end{align*}
All summands are irreducible, and are pairwise non-isomorphic.
In particular, $I(\chi)$ is not multiplicity-free.
\end{theorem}

\begin{proof}
By Lemma~\ref{lem:compare}(1), $\chi^{\delta_3,\psi}=\chi^{\prime,\psi}$ since $\chi_1=\chi_2$, and $\chi^{\prime,\psi}\neq\chi^{*}_{\psi}$ since $\chi_1\neq\chi_3$. 
The decomposition follows from Lemma~\ref{lem:P21} and the fact that $D(\chi)=10$.
\end{proof}

\begin{remark}\label{rem:D10-other}
Similar statements hold in the other two cases with $D(\chi)=10$, namely $\chi_2=\chi_3\neq\chi_1$ and $\chi_1=\chi_3\neq\chi_2$ (with $\chi_1\sim\chi_2\sim\chi_3$), and we omit the details.
In each case $I(\chi)$ has $7$ distinct irreducible constituents, and exactly one of them occurs with multiplicity two: it is $\Ind_{B_{\delta_3}}^G\chi^{\delta_3,\psi}$ if $\chi_2=\chi_3$, since then $\chi^{*}_{\psi}=\chi^{\delta_3,\psi}$, and $\Ind_{B'}^G\chi^{\prime,\psi}$ if $\chi_1=\chi_3$, since then $\chi^{*}_{\psi}=\chi^{\prime,\psi}$.
\end{remark}

The proof of the following lemma is similar to that of Lemma~\ref{lem:P21}, and we omit it.

\begin{lemma}\label{lem:alleq}
Let $\chi_1=\chi_2=\chi_3$. Then the character $\ct$ of $B_1$ extends to the character $\ct=\chi_1\circ\det$ of $G$, and
\[
\Ind_{B_1}^G\ct\;\cong\;\ct\;\oplus\;2\,(\ct\otimes\pi_1)\;\oplus\;(\ct\otimes\St_3),
\]
where $\St_3$ is the Steinberg representation of $\GL_3(\Fq)$ inflated to $G$, and $\pi_1$ is the irreducible representation of $\GL_3(\Fq)$ of dimension $q^2+q$ occurring in $\Ind_{\overline{B}}^{\GL_3(\Fq)}\one$, inflated to $G$.
\end{lemma}

\begin{theorem}\label{thm:D18}
Let $\chi_1=\chi_2=\chi_3$. Then $\chi^{\delta_3,\psi}=\chi^{\prime,\psi}=\chi^{*}_{\psi}$, and
\begin{align*}
I(\chi)\cong\ct\oplus2\,(\ct\otimes\pi_1)\oplus(\ct\otimes\mathrm{St}_3)\oplus3\cdot\Ind_{B_{\delta_3}}^G\chi^{\delta_3,\psi}
\oplus\Ind_{B''}^G\chi^{\prime\prime,\psi}\oplus\Ind_{B''}^G\chi^{**}_{\psi}\oplus\Ind_{B^{(2)}_{\delta_1}}^G\chi^{(2,\psi)}_{\delta_1},
\end{align*}
where $\ct=\chi_1\circ\det$ is as in Lemma~\ref{lem:alleq}.
All summands are irreducible and pairwise non-isomorphic.
In particular $I(\chi)$ is not multiplicity-free.
\end{theorem}

\begin{proof}
By Lemma~\ref{lem:compare}(1) the three characters of $B'$ coincide, and the decomposition follows from Lemma~\ref{lem:alleq}.
The constituents of $\Ind_{B_1}^G\ct$ are $\ct$, $\ct\otimes\pi_1$ and $\ct\otimes\mathrm{St}_3$, which are irreducible and pairwise non-isomorphic, and contribute $1+2^2+1+3^2+1+1+1=18$ to $D(\chi)$ (Lemma~5.10).
This proves the theorem.
\end{proof}

\subsection{The regular constituent of \texorpdfstring{$I(\chi)$}{I(chi)}}\label{subsec:regular}
Let $G=\GL_3(\cO_2)$, $K=I_3+\vp M_3(\cO_2)$, and $\psi_\beta(I_3+\vp Y)=\psi(\vp\tr(\beta Y))$ for $\beta\in M_3(\Fq)$. 
For $g\in G$ with image $\bar g\in\GL_3(\Fq)$ we have $\psi_\beta(g^{-1}kg)=\psi_{\bar g\beta\bar g^{-1}}(k)$. Let $\beta\in M_3(\Fq)$ be regular, let $\tilde\beta\in M_3(\cO_2)$ be a lift, and put $A=\cO_2[\tilde\beta]$. 
The centralizer of $\beta$ in $\GL_3(\Fq)$ is $\Fq[\beta]^\times$, so the stabilizer of $\psi_\beta$ in $G$ is $I(\psi_\beta)=A^\times K$, and $I(\psi_\beta)/K\cong\Fq[\beta]^\times$ is abelian. Using the construction of regular representations of $\GL_3(\cO_2)$ given in~\cite{Hill95,KOS18}, there exists $\tilde\psi_\beta$ an extension of $I(\psi_\beta)$ such that $\rho:=\Ind_{I(\psi_\beta)}^G\tilde\psi_\beta$ is an irreducible regular representation of $\GL_3(\cO_2)$ with dimension $\dim\rho=[G:I(\psi_\beta)]=[\GL_3(\Fq):C_{M_3(\Fq)}(\beta)^\times]$, where $C_{M_3(\Fq)}(\beta)$ is the centraliser of $\beta$ in $M_3(\Fq)$.

Write $\chi_i(1+\vp x)=\psi(\vp c_ix)$ with $c_i\in\Fq$, and put
\[
\beta=\begin{pmatrix}c_1&1&0\\ 0&c_2&1\\ 0&0&c_3\end{pmatrix},\qquad\tilde\beta=\begin{pmatrix}\tilde c_1&1&0\\ 0&\tilde c_2&1\\ 0&0&\tilde c_3\end{pmatrix},
\]
where $\tilde c_i\in\cO_2$ is a lift $c_i$ for $i=1,2,3$. Then $\beta$ is regular element of $M_3(\Fq)$ and put $A=\cO_2[\tilde\beta]$. Then $I(\psi_\beta)=A^\times\cdot K$, where $I(\psi_\beta)/K$ is the centraliser of $\beta\in\GL_3(\Fq)$. Define
\[
\lambda_\chi(ak):=\chi(a)\,\psi_\beta(k)\qquad(a\in A^\times,\ k\in K).
\]

\begin{theorem}\label{thm:regular}
$\lambda_\chi$ is a well-defined character of $I(\psi_\beta)=A^\times K$, and $\rho_\chi:=\Ind_{I(\psi_\beta)}^G\lambda_\chi$ is the unique irreducible regular constituent of $I(\chi)$. It occurs in $I(\chi)$ with multiplicity one.
\end{theorem}

\begin{proof}
We first prove that $\lambda_\chi$ is a character.
The group $A^\times$ normalizes $K$ and fixes the character $\psi_\beta$, and $\chi|_{A^\times}$ is a character because $A^\times\subset B$.
It remains to check that $\chi$ and $\psi_\beta$ agree on $A^\times\cap K$. 
For $a=I_3+\vp Y\in A^\times\cap K$ with $Y=(Y_{ij})$ upper triangular, so $\tr(\tilde{\beta} Y)=\sum_i\tilde{c}_iY_{ii}$.
Hence
\[
\psi_\beta(a)=\psi\Bigl(\vp\sum_i\tilde{c}_iY_{ii}\Bigr)=\prod_i\chi_i(1+\vp Y_{ii})=\chi(a).
\]
Now, we prove that $\rho_\chi$ is a sub representation of $I(\chi)$. Let $\lambda$ be any character of $I(\psi_\beta)$ with $\lambda|_K=\psi_\beta$. By Mackey theory,
\begin{equation}\label{eq--1}
\Hom_G\bigl(\Ind_{I(\psi_\beta)}^G\lambda,\ I(\chi)\bigr)\cong\bigoplus_{x\in[I(\psi_\beta)\backslash G/B]}\Hom_{I(\psi_\beta)\cap xBx^{-1}}\bigl(\lambda,\chi^x\bigr).
\end{equation}
Suppose for some $x\in[I(\psi_\beta)\backslash G/B]$, $\Hom_{I(\psi_\beta)\cap xBx^{-1}}\bigl(\lambda,\chi^x\bigr)\neq0$. Since $K$ is normal in $G$, we have $K\cap xBx^{-1}=x(K\cap B)x^{-1}$. Therefore,
\[
\Hom_{I(\psi_\beta)\cap xBx^{-1}}\bigl(\lambda,\chi^x\bigr)\neq0\implies\Hom_{K\cap xBx^{-1}}\bigl(\lambda,\chi^x\bigr)\neq0.
\]
We get $\psi_\beta(x(I_3+\vp Y)x^{-1})=\chi(I_3+\vp Y)$ for all upper triangular $Y$, i.e.
\[
\psi\bigl(\vp\tr(x^{-1}\tilde{\beta} xY)\bigr)=\psi\Bigl(\vp\sum_i\tilde{c}_iY_{ii}\Bigr).
\]
Since $\psi$ is a non-trivial character of $\mathcal{O}_2$ such that $\psi$ is non-trivial on $\varpi \mathcal{O}_2$ and using the fact that $\varpi\mathcal{O}_2\cong \F_q$, we treat the character $\psi$ as a character of $\F_q$ and get 
\[
\psi\bigl(\tr(\bar{x}^{-1}\beta\bar{x}\bar{Y})\bigr)=\psi\Bigl(\sum_ic_i\bar{Y}_{ii}\Bigr).
\]
Hence $\bar{x}^{-1}\beta \bar{x}$ is upper triangular with diagonal $(c_1,c_2,c_3)$.
In fact, it is simple to see that $\bar x\in B(\Fq)$, i.e. $x\in KB\subset I(\psi_\beta)B$.
Hence \(\Hom_{I(\psi_\beta)\cap xBx^{-1}}\bigl(\lambda,\chi^x\bigr)\neq 0\implies  x=I_3\).

For $x=I_3$, since $A^\times\subset B$ we have $I(\psi_\beta)\cap B=A^\times(K\cap B)$. 
On $K\cap B$ the characters $\lambda$ and $\chi$ agree by (i), so $\Hom_{I(\psi_\beta)\cap B}(\lambda,\chi)\neq0$ if and only if $\lambda|_{A^\times}=\chi|_{A^\times}$, i.e., $\lambda=\lambda_\chi$. 
Therefore
\[
\dim\Hom_G\bigl(\Ind_{I(\psi_\beta)}^G\lambda,\ I(\chi)\bigr)=\begin{cases}1,&\lambda=\lambda_\chi,\\ 0,&\lambda\neq\lambda_\chi.\end{cases}
\]

Now, we prove that $I(\chi)$ does not contain other regular representation. 
By Mackey theory for the normal subgroup $K$,
\[
I(\chi)|_K\cong\bigoplus_{x\in G/KB}\Ind_{K\cap xBx^{-1}}^K\chi^x.
\]
The above computation shows that for $\gamma\in M_3(\Fq)$, a character $\psi_\gamma$ occurs in $I(\chi)|_K$ only if $\bar x^{-1}\gamma\bar x$ is upper triangular with diagonal $(c_1,c_2,c_3)$ for some $x$. 
Hence every regular constituent of $I(\chi)$ lies over the $G$-orbit of $\psi_\beta$, so it is of the form $\rho_\chi$ and multiplicity of $\rho_\chi$ in $I(\chi)$ is one.
\end{proof}

\begin{corollary}\label{cor:regdim}
The dimension of the regular representation $\rho_\chi$ in $I(\chi)$ is given by
\[
\dim\rho_\chi=[\GL_3(\Fq):C_{M_3(\Fq)}(\beta)^\times]=\begin{cases}
q^3(q+1)(q^2+q+1),&\text{if the }c_i\text{ are pairwise distinct},\\
q^2(q+1)(q^3-1),&\text{if exactly two of the }c_i\text{ are equal},\\
q(q^2-1)(q^3-1),&\text{if }c_1=c_2=c_3.
\end{cases}
\]
Equivalently, these three cases correspond, respectively, to the $\chi_i$ being pairwise $\not\sim$, exactly two of them being equivalent under $\sim$, and $\chi_1\sim\chi_2\sim\chi_3$.
\end{corollary}

\section{Acknowledgements}
The author thanks Dipendra Prasad for carefully reading the manuscript and for his valuable comments.
The author had used the argument of Lemma~2.4 in earlier work without stating it precisely; its statement in the present form was formulated with the help of AI tools, which were also used to improve the language and clarity of the manuscript.

\end{document}